\documentclass[a4paper,12pt,final,reqno]{amsart}

\usepackage{amsmath,amssymb,amsthm}

\numberwithin{equation}{section}
\allowdisplaybreaks

\newtheorem{thm}{Theorem}[section]
\newtheorem{prp}{Proposition}[section]
\newtheorem{lem}{Lemma}[section]
\newtheorem{rmk}{Remark}[section]
\newtheorem{dfn}{Definition}[section]
\newtheorem{cor}{Corollary}[section]
\newtheorem{con}{Conjecture}[section]

\newcommand{\bbC}{\mathbb{C}}

\newcommand{\bbZ}{\mathbb{Z}}

\title{Some transformation formulas associated with  Askey-Wilson polynomials
and Lassalle's formulas for Macdonald-Koornwinder polynomials}
\author{A.~Hoshino, M.~Noumi and J.~Shiraishi}
\address{AH: Kagawa National College of Technology, 
355 Chokushi-cho, Takamatsu, Kagawa761-8058, Japan}
\email{hoshino@t.kagawa-nct.ac.jp}
\address{MN: Department of Mathematics, Kobe University, Rokko, Kobe 657-8501, Japan}
\email{noumi@math.kobe-u.ac.jp}
\address{JS: Graduate School of Mathematical Sciences, University of Tokyo, Komaba, Tokyo 153-8914, Japan}
\email{shiraish@ms.u-tokyo.ac.jp}

\begin{document}

\begin{abstract}
We present a fourfold series expansion representing the Askey-Wilson polynomials. 
To obtain the result,  
a sequential use is made of several summation and transformation formulas for the basic hypergeometric series,
including the Verma's $q$-extension of the Field and Wimp expansion, 
Andrews' terminating $q$-analogue of Watson's ${}_3F_2$ sum, 
Singh's quadratic transformation. As an application, we present an explicit formula for 
the Koornwinder polynomial of type $BC_n$ ($n\in \bbZ_{>0}$) with one row diagram. When the parameters are specialized, 
we recover Lassalle's formula for Macdonald polynomials of type $B_n$, $C_n$ and $D_n$ with one row diagram,
thereby proving his conjectures.
\end{abstract}

\maketitle

\section{Introduction}
Let $a,b,c,d,q \in \bbC $ be parameters with the condition $|q|<1$.
Let $D$ denote the Askey-Wilson $q$-difference operator \cite{AW}
\begin{align}
D&={(1-ax)(1-bx)(1-cx)(1-dx)\over (1-x^2)(1-q x^2)} \left(T_{q,x}^{+1}-1\right)\\
&+{(1-a/x)(1-b/x)(1-c/x)(1-d/x)\over (1-1/x^2)(1-q /x^2)} \left(T_{q,x}^{-1}-1\right)\nonumber,
\end{align}
where the $q$-shift operators are defined by  $T_{q,x}^{\pm1}f(x)=f(q^{\pm 1}x)$.
Recall the fundamental facts about the Askey-Wilosn polynomial $p_n(x;a,b,c,d|q)$ ($n\in \bbZ_{\geq 0}$). 
It is a symmetric Laurent polynomial in $x$ and characterized by the two conditions:
(i) $p_n(x)$ has the highest degree $n$, (2) $p_n(x)$ is an eigenfunction of the operator $D$.
Askey-Wilson's celebrated formula reads \cite{AW}
\begin{align}
&p_n(x)=a^{-n}(ab,ac,ad;q)_n\,\,
{}_4\phi_3\left[ {q^{-n}, abcdq^{n-1},ax,a/x\atop ab,ac,ad};q,q \right],\label{AWpn}\\
&D p_n(x)= \left(q^{-n}+qbcdq^{n-1}-1-abcdq^{-1}\right)p_n(x). \label{Dpn}
\end{align}
Here and hereafter we use the standard notations (see \cite{GR} for more details)
\begin{align}
&(a_1,a_2,\ldots,a_k;q)_n= \prod_{j=1}^k \prod_{i=0}^{n-1}(1-q^i a_j),\\
&{}_{r+1}\phi_r\left[{a_1,a_2,\ldots,a_{r+1}\atop b_1,b_2,\ldots,b_{r}};q,x \right]
=\sum_{m\geq 0} {(a_1,a_2,\ldots,a_{r+1};q)_m \over (q,b_1,b_2,\ldots,b_{r};q)_m} x^m.
\end{align}

Let $s\in \bbC$ be a parameter. Introduce $\lambda$ satisfying $s=q^{-\lambda}$.
Then we have $T_{q,x} x^{-\lambda} =s x^{-\lambda}$.
Let $f(x;s)=f(x;s|a,b,c,d|q)$ be a formal series in $x$
\begin{align}
&f(x;s)=x^{-\lambda} \sum_{n\geq 0} c_n x^{n},\quad c_0\neq 0,\label{cn}
\end{align}
satisfying the $q$-difference equation
\begin{align}
&D f(x;s)=\left(s+{abcd\over qs}-1-{abcd\over q}\right)f(x;s). \label{Df}
\end{align}
With the normalization $c_0=1$, (\ref{Df})
determines the coefficients $c_n=c_n(s|a,b,c,d|q)$ uniquely as rational functions in $a,b,c,d,q$ and $s$.

By using (\ref{AWpn}), we can easily find an explicit formula for $f(x;s)$. (See Section \ref{6-phi-5}, Theorem \ref{6phi5}.)
\begin{thm}\label{thm-1}
We have
\begin{align}
&f(x;s)= x^{-\lambda}{(ax;q)_\infty\over (qx/a;q)_\infty}
\sum_{n\geq 0}
{(qs^2/a^2;q)_n \over (q;q)_n} (a x/s)^n \label{6phi5eqfintro}\\
&\times
{}_6\phi_5
\left[
{q^{-n},q^{n+1}s^2/a^2,s,qs/ab,qs/ac,qs/ad 
\atop
q^{2}s^2/abcd,q^{1/2}s/a,-q^{1/2}s/a,qs/a,-qs/a
};q,q
\right]. \nonumber
\end{align}
\end{thm}

\begin{rmk}
When $\lambda=n\in \bbZ_{\geq 0}$, the series $f(x,q^{-n})$ must be 
proportional to the Askey-Wilson polynomial $p_n(x)$, in particular, indicating the termination of the series.
Note that, 
however, such termination can not easily be seen from the expression (\ref{6phi5eqfintro}).
This is one of the reason that we seek another expression below. 
\end{rmk}

Studying the series $f(x;s)$, one finds
that several interesting techniques are involved
including the Verma's $q$-extension of the Field and Wimp expansion, 
Andrews' terminating $q$-analogue of Watson's ${}_3F_2$ sum, 
Singh's quadratic transformation \cite{GR}. (See Section \ref{PROOF}.)

\begin{dfn}
 Set $c_e(k,l;s)=c_e(k,l;s|a,c|q)$ and $c_e(k,l;s)=c_e(k,l;s|a,b,c,d|q)$ by
\begin{align}
c_e(k,l;s)&=
{(a^2;q^2)_k (q^{4l}s^2;q^2)_k \over (q^2;q^2)_k (q^{4l} q^2 s^2/a^2;q^2)_k}(q^2/a^2)^k\label{ce}\\ 
&\times
{(c^2/q;q^2)_l (s^2/a^2;q^2)_l \over 
(q^2;q^2)_l (q^3 s^2/a^2c^2;q^2)_l}
{(s;q)_{2l} (q^2s^2/a^4;q^2)_{2l} \over(qs/a^2;q)_{2l}  (s^2/a^2;q^2)_{2l}}(q^2/c^2)^l,\nonumber\\
c_o(m,n;s)&=
{(-b/a;q)_m(s;q)_m(qs/cd;q)_m(qs^2/a^2c^2;q)_m\over 
(q;q)_m(q^2s^2/abcd;q)_m(qs^2/a^2c^2;q^2)_m}(q/b)^m\label{co}\\
&\times 
{(-d/c;q)_n(q^m s;q)_n(qs/ab;q)_n(-q^m qs/ac;q)_n (q^m qs^2/a^2c^2;q)_n\over 
(q;q)_n(q^m q^2s^2/abcd;q)_n (-qs/ac;q)_n(q^{2m}qs^2/a^2c^2;q^2)_n}(q/d)^n.\nonumber
\end{align}
\end{dfn}
\begin{rmk}
(1) The `even generators' $c_e(k,l;s)$ are basically composed in terms of
the $q$-shifted factorials with the base $q^2$, and the `odd ones' $c_o(k,l;s)$ are with the base $q$.
(2) The $c_e(k,l;s)$ does not depend on $b$ and $d$.
(3) The $c_o(m,n;s)$ is recast as
\begin{align}
c_o(m,n;s)&=
{(-b/a;q)_m(qs/cd;q)_m\over 
(q;q)_m(-qs/ac;q)_m}
{(s;q)_{m+n}(-qs/ac;q)_{m+n}(qs^2/a^2c^2;q)_{m+n}\over 
(q^2s^2/abcd;q)_{m+n}(q^{1/2}s/ac ;q)_{m+n}(-q^{1/2}s/ac ;q)_{m+n}}
(q/b)^m \nonumber\\
&\times 
{(-d/c;q)_n (qs/ab;q)_n\over 
(q;q)_n (-qs/ac;q)_n}(q/d)^n.
\end{align}
\end{rmk}
Now we state our main result in the present paper. (See Section \ref{PROOF}, Theorems \ref{thm:special} 
and \ref{thm:general}.)
\begin{thm}\label{main}
Let $s\in \bbC$ be generic. 
We have
\begin{align}
f(x;s)=x^{-\lambda} \sum_{k,l,m,n\geq 0}  c_e(k,l;q^{m+n}s) c_o(m,n;s)x^{2k+2l+m+n}.\label{fourfold}
\end{align}
\end{thm}

This gives us a fourfold summation formula for the Askey-Wilson polynomial.
\begin{thm}
Let $\lambda\in \bbZ_{\geq 0}$. 
We have the formula for the Askey-Wilson polynomial $p_\lambda(x)$ representing 
as a sum of monomials in $x$ with factorized coefficients
\begin{align}
p_\lambda(x)=(abcdq^{\lambda-1};q)_\lambda \,\,
\sum_{(k,l,m,n)\in {\mathcal P}_\lambda} 
 c_e(k,l;q^{m+n-\lambda}) c_o(m,n;q^{-\lambda})x^{-\lambda+2k+2l+m+n},
\end{align}
where ${\mathcal P}_\lambda\subset (\bbZ_{\geq 0})^4$ denotes the 
finite set of points in the polyhedron defined by the set of inequalities 
\begin{align}
0\leq m\leq \lambda,\quad  0\leq n \leq \lambda-m,\quad 
0\leq 2l\leq \lambda -m-n,\quad 0\leq k\leq \lambda-2l-m-n.
\end{align}
\end{thm}

The even part of the series can be transformed 
by using Singh's formula or the $q$-analogue of Bailey's formula \cite{GR}. 
(See Propositions \ref{prp:type-2} and 
\ref{prp:type-3}.) 
\begin{prp}\label{PRP-1}
We have the two bibasic representations with bases $q$ and $q^2$:
\begin{align}
&\sum_{k,l\geq 0} c_e(k,l;s|a,c|q) x^{2k+2l}\\
&=
\sum_{k,l\geq 0} 
{(qa^2/c^2;q^2)_k(q^{2l}s^2;q^2)_k\over 
(q^2;q^2)_k(q^{2l} q^3s^2/a^2c^2;q^2)_k} (q^2x^2/a^2)^k 
{(c^2/q;q)_l(s;q)_l (q^2s^2/a^4;q^2)_l\over 
(q;q)_l(qs/a^2;q)_l (q^3 s^2/a^2c^2;q^2)_l} (q^2x^2/c^2)^l\nonumber \\
&=
\sum_{k,l\geq 0} 
{(q a^2/c^2;q^2)_k (q^3 s/c^2;q^2)_k (q^2 s^2/c^4;q^2)_k \over 
(q^2;q^2)_k(q s/c^2;q^2)_k(q^3 s^2/a^2 c^2;q^2)_k}(q^2x^2/a^2)^k 
{(c^2/q;q)_l (s;q)_{2k+l} \over (q;q)_l (q^2 s/c^2;q)_{2k+l} }(q^2x^2/c^2)^l.  \nonumber
\end{align}
\end{prp}

As an application, 
we present an explicit formula for Koornwinder  polynomials \cite{Ko} with one row
diagram. It is derived from Theorem \ref{main} and Proposition \ref{PRP-1} by using the kernel function 
which intertwines the action of the Koornwinder operators of type $BC_n$ and $BC_1$ \cite{KNS}. 
 (See Section \ref{Koornwinder} below).   
Let $n\in \bbZ_{>0}$ and $x=(x_1,\ldots,x_n)$ be a set of variables. 
Let $P_{(r)}(x|a,b,c,d|q,t)$ be the Koornwinder polynomial with one row diagram $(r)$ ($r\in \bbZ_{\geq 0}$).
 Set 
\begin{align}
g_r(x|a,b,c,d|q,t)={(t;q)_r\over (q;q)_r} P_{(r)}(x|a,b,c,d|q,t),
\end{align}
for simplicity of display.
\begin{dfn} \label{G_r}
Define  the symmetric Laurent polynomial $G_r(x;q,t)$ by 
\begin{align}
\prod_{i=1}^n {(tux_i;q)_\infty  \over(ux_i;q)_\infty  }
{(tu/x_i;q)_\infty  \over(u/x_i;q)_\infty  }=\sum_{r\geq 0} G_r(x;q,t) u^r.
\end{align}
\end{dfn}

\begin{thm}\label{thm:g_r}
We have
\begin{align}
&g_r(x|a,-a,c,-c|q,t)\\
&=\sum_{k,l\geq 0\atop 2k+2l\leq r}
G_{r-2k-2l}(x;q,t)
{(q a^2/c^2;q^2)_k (q^{3-r} t^{1-n} /c^2;q^2)_k (q^{2-2r} t^{2-2n}/c^4;q^2)_k \over 
(q^2;q^2)_k(q^{1-r} t^{1-n}/c^2;q^2)_k(q^{3-2r} t^{2-2n}/a^2 c^2;q^2)_k}(t^2/a^2)^k \nonumber \\
&\quad \quad \times
{(c^2/qt;q)_l (q^{-r}t^{-n}/q;q)_{2k+l} \over (q;q)_l (q^{2-r} t^{1-n}/c^2;q)_{2k+l} }
{1-q^{-r+2k+2l}t^{-n}\over 1-q^{-r}t^{-n}}
(t^2/c^2)^l ,\nonumber \\[4mm]
&g_r(x|a,b,c,d|q,t)\label{g_rabcd}\\
&=\sum_{i,j\geq 0\atop i+j\leq r}
g_{r-i-j}(x|a,-a,c,-c|q,t) 
{(-b/a;q)_i(q^{1-r}t^{1-n}/cd;q)_i\over 
(q;q)_i(-q^{1-r}t^{1-n}/ac;q)_i} \nonumber\\
&\quad\quad \times
{(q^{1-r}t^{-n};q)_{i+j}(-q^{1-r}t^{1-n}/ac;q)_{i+j}(r^{1-2r}t^{2-2n}/a^2c^2;q)_{i+j}\over 
(q^{2-2r}t^{2-2n}/abcd;q)_{i+j}(q^{1/2-r}t^{1-n}/ac ;q)_{i+j}(-q^{1/2-r}t^{1-n}/ac ;q)_{i+j}}
(t/b)^i \nonumber\\
&\quad\quad\times 
{(-d/c;q)_j (q^{1-r}t^{1-n}/ab;q)_j\over 
(q;q)_j (-q^{1-r}t^{1-n}/ac;q)_n}(t/d)^j. \nonumber
\end{align}
\end{thm}

\begin{cor}
By specializing the parameters in (\ref{g_rabcd}), we recover Lassalle's formulas for 
the Macdonald polynomials of type $B$, $C$ and $D$,
thereby proving his Conjectures 1, 3 and 4 in \cite{Lass}. 
\end{cor}

Note that for Macdonald polynomials of type $B$, $C$ and $D$, it is convenient to use
the following simplified version of the series $f(x;s)$. (See Propositions \ref{simplify-1} and \ref{simplify-2}.)
\begin{prp}\label{PRP-2}
We have
\begin{align}
&f(x;s|-a,b,-q^{1/2} a,q^{1/2} b|q)=x^{-\lambda}
\sum_{l,m\geq 0}
{(a^2;q)_m(q^l s;q)_m\over (q;q)_m(q^{l+1}s/a^2;q)_m} (q x^2/a^2)^m\\
&\quad \quad \times
{(b/a;q^{1/2})_l(s/a^2;q^{1/2})_l  (s;q)_l\over
 (q^{1/2};q^{1/2})_l (q^{1/2}s/ab;q^{1/2})_l (s/a^2;q)_l }(q^{1/2}x/b)^l, \nonumber\\
&f(x;s|-a,b,-q^{1/2} a,q^{1/2} a|q)=x^{-\lambda}
\sum_{l,m\geq 0}
{(a^2;q)_m(q^l s;q)_m\over (q;q)_m(q^{l+1}s/a^2;q)_m} (q x^2/a^2)^m\\
&\quad \quad \times
{(b/a;q)_l(s^2/a^4;q)_l  (s;q)_l\over
 (q;q)_l (qs^2/a^3b;q)_l (s/a^2;q)_l }(qx/b)^l. \nonumber
\end{align}
\end{prp}

The structure of this paper is as follows. 
In Section \ref{6-phi-5}, we give a proof of Theorem \ref{thm-1} (See Corollary \ref{cor-f}.).
In Section \ref{PROOF}, Theorem \ref{main} is proved in two steps, 
first for the special parameters $(a,b,c,d)=(a,-a,c,-c)$ (\S\S \ref{special}, Theorem \ref{thm:special}), 
and next for the general parameters $(a,b,c,d)$ (\S\S \ref{general}, Theorem \ref{thm:general}).
In Section \ref{BIBASIC}, we prove Propositions \ref{PRP-1} and \ref{PRP-2}
using some bibasic transformation formulas. Section \ref{Koornwinder} is devoted 
to the proof of Theorem \ref{thm:g_r}.
%
In \S\S \ref{Mac-C} and \S\S \ref{Mac-BD} the explicit formulas are derived 
for the Macdonald polynomials with one row for type $C_n$, $B_n$ and $D_n$,
thereby proving Lassalle's conjectures \cite[p.8, Conjecture 1, p.10, Conjecture 3, p.11, Conjecture 4]{Lass}. 
In Appendix, some basic facts are recalled concerning  the kernel function associated with 
Koornwinder's difference operator \cite{KNS}.
In \S\S \ref{reproduction} Theorem \ref{thm:reproduction},
a reproduction formula for the Koornwinder polynomials is given 
in terns of the kernel function.
Some notations 
for the Macdonald polynomials of type $C_n$, $B_n$ and $D_n$ are briefly given
in \S\S \ref{notation-C} and \S\S \ref{notation-BD}. In \S\S \ref{con-B2},
we present a conjecture about the Macdonald polynomial of $B_2$.  
\bigskip

\noindent
{\bf 
Acknowledgments.}
We are grateful to M.~Lassalle and V.~Pasquier for stimulating discussions.
Research of J.S. is supported by the Grant-in-Aid for Scientific
Research C-24540206.

\section{Proof of Theorem \ref{thm-1}}\label{6-phi-5}
For simplicity of display, we introduce a notation.
\begin{dfn}
Set 
\begin{align}
\Psi(x;s|a,b,c,d|q)&=
 {(ax;q)_\infty\over (qx/a;q)_\infty}
\sum_{n\geq 0}
{(qs^2/a^2;q)_n \over (q;q)_n} (a x/s)^n \label{6phi5eqf}\\
&\times
{}_6\phi_5
\left[
{q^{-n},q^{n+1}s^2/a^2,s,qs/ab,qs/ac,qs/ad 
\atop
q^{2}s^2/abcd,q^{1/2}s/a,-q^{1/2}s/a,qs/a,-qs/a
};q,q
\right]. \nonumber
\end{align}
\end{dfn}

\begin{thm}\label{6phi5}
We have
\begin{align}
{}_4\phi_3\left[ {q^{-m}, abcdq^{m-1},ax,a/x\atop ab,ac,ad};q,q \right]
=a^m x^{-m}
{(abcdq^{m-1};q)_m \over (ab,ac,ad;q)_m}\Psi(x;q^{-m}|a,b,c,d|q).\label{6phi5eq}
\end{align}
\end{thm}
This gives an explicit formula for the infinite series $f(x;s)$ for generic $s\in \bbC$.
\begin{cor}\label{cor-f}
For generic $s\in \bbC$, $f(x;s)=x^{-\lambda} \Psi(x;s|a,b,c,d|q)$ satisfies (\ref{Df}). Hence Theorem \ref{thm-1} holds.
\end{cor}

\proof
The coefficients $c_n$'s in (\ref{cn}) are clearly rational functions in $s$. 
After clearing the denominator $(1-x^2)(1-qx^2)(1-q^{-1}x^2)$ in
 (\ref{Df}), we have a set of linear relations for $c_n$'s with coefficients 
being polynomials in $s$. Hence by substituting the $c_n$'s
to these linear relations and clearing the denominators, we have a set of polynomials equations in $s$. 
Theorem \ref{6phi5} and (\ref{Dpn}) mean 
that these polynomials are zero for infinitely many points $s=q^{-n}$ ($n=0,1,2,\ldots$).
Hence all of such polynomials are identically zero, indicating that (\ref{Df}) holds for generic $s\in \bbC$.
\qed
\medskip

\noindent
{\it Proof of Theorem \ref{6phi5}.}
We need to expand the ${}_4\phi_3$ series in the form of $x^{-m}$ times a power series in $x$.
Therefore our starting point should be
\begin{align}
\mbox{LHS of (\ref{6phi5eq})}
=
\sum_{k\geq 0}
{(q^{-m},abcdq^{m-1},ax,a/x;q)_{m-k}\over (q,ab,ac,ad;q)_{m-k}}q^{m-k}.\label{start} 
\end{align}
By using the $q$-binomial formula \cite[p.7, (1.3.2)]{GR} we have
\begin{align}
&(a x;q)_{m-k}(a /x;q)_{m-k}
={(ax;q)_\infty\over (qx/a;q)_\infty} {(q^{-m+k+1}x/a;q)_\infty\over (q^{n-k}ax;q)_\infty}(-a/x)^{m-k}q^{(m-k)(m-k-1)/2} 
\nonumber\\
&={(ax;q)_\infty\over (qx/a;q)_\infty} (-a/x)^{m-k}q^{(m-k)(m-k-1)/2} 
\sum_{l\geq 0} {(q^{-2m+2k+1}/a^2;q)_l\over (q;q)_l} (q^{m-k}ax)^l.
\end{align} 
Then simplifying the factors we have
\begin{align}
&\mbox{RHS of (\ref{start})}\nonumber \\
&=a^m x^{-m} {(abcd q^{m-1};q)_{m} \over (ab,ac,ad;q)_{m}}{(ax;q)_\infty\over (qx/a;q)_\infty} 
\sum_{k\geq 0}
\sum_{l\geq 0}
{(q^{-m},q^{-m+1}/ab,q^{-m+1}/ac,q^{-m+1}/ad;q)_k \over (q,q^{-2m+2}/abcd;q)_k} \nonumber\\
&\quad \quad \times
{(q^{-2m+2k+1}/a^2;q)_l\over (q;q)_l} 
(-1)^k a^{k+l}x^{k+l} q^{(m-k)(k+l)+{k(k+1)\over 2}}\\
&=a^m x^{-m} {(abcd q^{m-1};q)_{m} \over (ab,ac,ad;q)_{m}}{(ax;q)_\infty\over (qx/a;q)_\infty} 
\sum_{n\geq 0}
{(q^{-2m+1}/a^2;q)_n \over (q;q)_n} (q^m a x)^n \nonumber\\
&\times
\sum_{k\geq 0}
{(q^{-m},q^{-m+1}/ab,q^{-m+1}/ac,q^{-m+1}/ad;q)_k \over (q,q^{-2m+2}/abcd;q)_k}
{(q^{-n},q^{-2m+n+1}/a^2;q)_k \over (q^{-2m+1}/a^2;q)_{2k}}q^k=\mbox{RHS of (\ref{6phi5eq})}.\nonumber
\end{align}
\qed

\section{Proof of Theorem \ref{main}}\label{PROOF}
We embark on the proof of Theorem \ref{main}. For clarity of display we need another notation.
\begin{dfn}
Set 
\begin{align}
\Phi(x;s|a,b,c,d|q)= \sum_{k,l,m,n\geq 0}  c_e(k,l;q^{m+n}s|a,c|q) c_o(m,n;s|a,b,c,d|q)x^{2k+2l+m+n}.
\end{align}
\end{dfn}
In view of Corollary \ref{cor-f},
we only need to show that $\Phi(x;s|a,b,c,d|q)=\Psi(x;s|a,b,c,d|q)$.
We shall divide our proof in two steps. 
First, we will consider the special case $b=-a,d=-c$ in 
\S\S  \, \ref{special}.  
Then the the general case will be treated in \S\S \, \ref{general}.

We remark that 
the special case $b=-a,d=-c$ (Theorem \ref{thm:special} below) 
constitute  the essential part of the proof of Theorem \ref{main}. 
Based on Theorem \ref{thm:special}, the general case can be treated easily (Theorem \ref{thm:general} below).

\subsection{Case $b=-a,d=-c$}\label{special}
For the sake of clarity we first write down $\Phi(x;s|a,-a,c,-c|q)$ and $\Psi(x;s|a,-a,c,-c|q)$ explicitly
\begin{align}
&\Phi(x;s|a,-a,c,-c|q)=\sum_{k,l\geq 0}c_e(k,l;s)x^{2k+2l},\\
&\Psi(x;s|a,-a,c,-c|q)
={(ax;q)_\infty\over (qx/a;q)_\infty}
\sum_{n\geq 0}
{(qs^2/a^2;q)_n \over (q;q)_n} (a x/s)^n \label{6phi5even} \\
&\quad \quad \times
{}_6\phi_5
\left[
{q^{-n},q^{n+1}s^2/a^2,s,-qs/a^2,qs/ac,-qs/ac 
\atop
q^{2}s^2/a^2c^2,q^{1/2}s/a,-q^{1/2}s/a,qs/a,-qs/a
};q,q
\right]. \nonumber
\end{align}

\begin{thm}\label{thm:special}
We have 
\begin{align}
\Phi(x;s|a,-a,c,-c|q)=\Psi(x;s|a,-a,c,-c|q).\label{6phi5even} 
\end{align}
\end{thm}

Recall  
Verma's $q$-extension of the Field and Wimp expansion  \cite[p.76, (3.7.9)]{GR}
\begin{align}
&{}_{r+t} \phi_{s+u}\left[{a_R,c_T \atop b_S,d_U};q,xw
\right]\label{Verma}\\
&=\sum_{j=0}^\infty 
{(c_T,e_K;q)_j \over (q,d_U,\gamma q^j;q)_j} x^j [(-1)^j q^{\left( j\atop 2\right)}]^{u+3-t-k}\nonumber \\
&\times 
{}_{t+k} \phi_{u+1}\left[{c_T q^j,e_K q^j  \atop \gamma q^{2j+1},d_U q^j};q,xq^{j(u+2-t-k)}
\right] 
{}_{r+2} \phi_{s+k}\left[{q^{-j},\gamma q^j,a_R\atop b_S,e_K};q,wq
\right].\nonumber
\end{align}
Here we have used the contracted notation $a_R$ for $a_1,\ldots,a_r$ by , etc.

Set the parameters in (\ref{Verma}) as
\begin{align}
&r=2,\quad s=2, \quad t=4, \quad  u=3,\quad k=1,\label{rstuk}\\
&w=1,\quad x=q,\\
&a_R=(qs/ac,-qs/ac),\qquad b_S=(q^{1/2}s/a,-q^{1/2}s/a),\nonumber\\
&c_T=(q^{-n},q^{n+1}s^2/a^2,s,-qs/a^2),\qquad d_U=(qs/a,-qs/a,q^2s^2/a^2c^2), \label{abcde}\\
&e_K=q^2s^2/a^2c^2,\qquad \gamma=s^2/a^2. \nonumber
\end{align}
Then the choice of the parameters (\ref{rstuk}) means that we have
an expansion of the form ${}_6\phi_5=\sum {}_5\phi_4\cdot {}_4\phi_3.$
Note, however, that  this ${}_5\phi_4$ series  degenerates to a ${}_4\phi_3$ series from the conditions (\ref{abcde}).

\begin{lem}\label{lem-1}
We have
\begin{align}
&{}_6\phi_5
\left[
{q^{-n},q^{n+1}s^2/a^2,s,-qs/a^2,qs/ac,-qs/ac 
\atop
q^{2}s^2/a^2c^2,q^{1/2}s/a,-q^{1/2}s/a,qs/a,-qs/a
};q,q
\right]\label{65-5443}\\
&=
\sum_{j\geq0}
{(q^{-n},q^{n+1}s^2/a^2,s,-qs/a^2,q^2s^2/a^2c^2;q)_j\over
(q,qs/a,-qs/a,q^2s^2/a^2c^2,q^j s^2/a^2;q)_j} (-1)^j q^{j+\left( j\atop 2\right)}\nonumber\\
&\times 
{}_4\phi_3
\left[
{q^{-n+j},q^{j+n+1}s^2/a^2,q^j s,-q^{j+1}s/a^2
\atop
q^{2j+1}s^2/a^2,q^{j+1}s/a,-q^{j+1}s/a
};q,q
\right] 
{}_4\phi_3
\left[
{q^{-j}q^j s^2/a^2,qs/ac,-qs/ac
\atop
q^{1/2}s/a,-q^{1/2}s/a,q^2s^2/a^2c^2
};q,q
\right].\nonumber
\end{align}
\end{lem}

We need to transform the two ${}_4\phi_3$ series on RHS of (\ref{65-5443}).
For the treatment of the last factor, 
we recall Andrews' terminating $q$-analogue of Watson's ${}_3F_2$ sum \cite[p.237, (II.17)]{GR}. 
Namely, we have
\begin{align}
{}_4\phi_3\left[ {q^{-n},aq^n,c,-c \atop (aq)^{1/2},-(aq)^{1/2},c^2};q,q\right]=
\left\{
\begin{array}{ll}
0,&\mbox{if $n$ is odd},\\[3mm]
\displaystyle
{c^n(q,aq/c^2;q^2)_{n/2}\over (aq,c^2 q;q^2)_{n/2}},& \mbox{if $n$ is even}.
\end{array}
\right.
\end{align}

\begin{lem} \label{lem-2}
From Andrews' formula, we have
\begin{align}
&{}_4\phi_3\left[ 
{q^{-j},q^j s^2/a^2,qs/ac,-qs/ac\atop
q^{1/2}s/a,-q^{1/2}s/a,q^2s^2/a^2c^2}
;q,q\right]\\
&=
\left\{
\begin{array}{ll}
0,&\mbox{if $j$ is odd},\\[3mm]
\displaystyle
\left(qs\over ac\right)^j
{(q,c^2/q;q^2)_{j/2}\over (qs^2/a^2,q^3s^2/a^2c^2;q^2)_{j/2}}& \mbox{if $j$ is even}.
\end{array}
\right. \nonumber
\end{align}
Hence the support of the summand in (\ref{65-5443}) is restricted to $j=0,2,4,\ldots$. 
\end{lem}
For the sake of simplicity, 
 we shall change our notation in (\ref{65-5443})
 as $j\rightarrow 2j$ and take summation over  $j=0,1,2,\ldots$.

Next, we turn to the other ${}_4\phi_3$ series.
Recall Singh's quadratic transformation \cite[p.89, (3.10.13)]{GR}, that is
\begin{align}
{}_4\phi_3
\left[
{a^2,b^2,c,d
\atop
abq^{1/2},-abq^{1/2},-cd
};q,q
\right]=
{}_4\phi_3
\left[
{a^2,b^2,c^2,d^2
\atop
a^2b^2q,-cd,-cdq
};q^2,q^2
\right],
\end{align}
provided the series terminates.

\begin{lem}\label{Singh}
For $n,j\in \bbZ_{\geq 0}$ satisfying $2j\leq n$,
we have the identity between the two terminating ${}_4\phi_3$ series
\begin{align}
&{}_4\phi_3
\left[
{q^{-n+2j},q^{2j+n+1}s^2/a^2,q^{2j}s,-q^{2j+1}s/a^2
\atop
q^{4j+1}s^2/a^2,q^{2j+1}s/a,-q^{2j+1}s/a
};q,q
\right]\label{Singheq}\\
&=
{(q/a^2;q)_{n-2j}\over (q^{4j+1}s^2/a^2;q)_{n-2j}} \left(q^{2j}s\right)^{n-2j}
{}_4\phi_3
\left[
{q^{-n+2j},q^{-n+2j+1},q^{4j}s^2,a^2
\atop
q^{4j+2}s^2/a^2,q^{-n+2j}a^2,q^{-n+2j+1}a^2
};q^2,q^2
\right].\nonumber
\end{align}
\end{lem}

\noindent
{\it Proof of Lemma \ref{Singh}.}
It is enough to show (\ref{Singheq}) for the infinite discrete points $q^{2j}s=q^{-k}$, $(k=0,1,2,\ldots)$.
Under this condition, we can apply Singh's formula \cite[p.89, (3.10.13)]{GR}, and we have
\begin{align}
&\mbox{LHS of (\ref{Singheq})}
=
{}_4\phi_3
\left[
{q^{-n+2j},q^{2j+n+1}s^2/a^2,q^{4j}s^2,q^{4j+2}s^2/a^4
\atop
q^{4j+2}s^2/a^2,q^{4j+1}s^2/a^2,q^{4j+2}s^2/a^2
};q^2,q^2
\right]\\
&
=
{(q^{2j+n+1}s^2/a^2,q^{2j+n+2}s^2/a^2;q^2)_k \over 
(q^{4j+1}s^2/a^2,q^{4j+2}s^2/a^2;q^2)_k } (q^{-n+2j})^k 
{}_4\phi_3
\left[
{q^{-n+2j},q^{-n+2j+1},q^{4j}s^2,a^2
\atop
q^{4j+2}s^2/a^2,q^{-n+2j}a^2,q^{-n+2j+1}a^2
};q^2,q^2
\right],\nonumber
\end{align}
where we used Sears' transformation \cite[p.41, (2.10.4)]{GR} in the second line. 
By noting $q^{2j}s=q^{-k}$, 
the factor in front of the ${}_4\phi_3$ series can be calculated as
 \begin{align}
(q^{-n+2j})^k
{(q^{2j+n+1}s^2/a^2;q)_k \over 
(q^{4j+1}s^2/a^2;q)_k } 
=(q^{2j}s)^{n-2j}{(q/a^2;q)_{n-2j}\over (q^{4j+1}s^2/a^2;q)_{n-2j}} .
\end{align}
\qed

Summarizing the results in Lemmas \ref{lem-1}, \ref{lem-2} and \ref{Singh}, 
we arrive at the following intermediate expression.
\begin{lem}\label{lem-4}
We have
\begin{align}
\mbox{RHS of (\ref{6phi5even})}
&=
{(ax;q)_\infty\over (qx/a;q)_\infty}
\sum_{n\geq 0}
\sum_{j= 0}^{\lfloor {n\over 2}\rfloor}
\sum_{m= 0}^{\lfloor {n-2j\over 2}\rfloor} A(n,j,m),
\end{align}
where 
\begin{align}
A(n,j,m)&={(q s^2/a^2;q)_n \over (q;q)_n} (a x/s)^n
{(q^{-n},q^{n+1}s^2/a^2,s,-qs/a^2;q)_{2j} \over 
(q,qs/a,-qs/a,q^{2j}s^2/a^2;q)_{2j}} q^{j(2j+1)} \nonumber\\
&\times 
{(q/a^2;q)_{n-2j} \over (q^{4j+1}s^2/a^2;q)_{n-2j}} (q^{2j }s)^{n-2j}\times
{(q,c^2/q;q^2)_j\over (q s^2/a^2,q^3s^2/a^2c^2;q^2)_j} \left(qs\over ac\right)^{2j}\\
&\times 
{(q^{-n+2j},q^{-n+2j+1},q^{4j}s^2,a^2;q^2)_{2m} .\nonumber
\over 
(q^2,q^{4j+2}s^2/a^2,q^{-n+2j}a^2,q^{-n+2j+1}a^2;q^2)_{2m}}q^{2m}.
\end{align}
Hence we can change the order of the summation as
\begin{align}
\sum_{n\geq 0}
\sum_{j= 0}^{\lfloor {n\over 2}\rfloor}
\sum_{m= 0}^{\lfloor {n-2j\over 2}\rfloor} A(n,j,m)=
\sum_{l\geq 0}\sum_{j\geq 0}\sum_{m\geq 0}A(l+2j+2m,j,m).
\end{align}

\end{lem}

\begin{lem}\label{lem-5}
We have
\begin{align}
&
{A(l+2j+2m,j,m)\over A(2j+2m,j,m)}=
{(q/a^2;q)_l\over (q;q)_l} (ax)^l,\\
&
A(2j+2m,j,m)=c_e(m,j;s)x^{2m+2j}.
\end{align}
\end{lem}
\proof Straightforward calculation. \qed

Now we are ready to present the final step.

\noindent
{\it Proof of Theorme \ref{thm:special}}
{} From Lemmas \ref{lem-4}, \ref{lem-5} and the $q$-binomial formula, we have
\begin{align}
\mbox{RHS of (\ref{6phi5even})}
&={(ax;q)_\infty\over (qx/a;q)_\infty} 
\sum_{l\geq 0}\sum_{j\geq 0}\sum_{m\geq 0}A(l+2j+2m,j,m)\\
&=\sum_{j\geq 0}\sum_{m\geq 0}A(2j+2m,j,m)=
\Phi(x;s|a,-a,c,-c|q)=
\mbox{LHS of (\ref{6phi5even})}.\nonumber
\end{align}
\qed

\subsection{General Case}\label{general}
Based on Theorem \ref{thm:special}, we treat the general case. 
\begin{thm}\label{thm:general}
We have
\begin{align}
\Psi(x;s|a,b,c,d|q)=\Phi(x;s|a,b,c,d|q),
\end{align}
which implies Theorem \ref{main}.
\end{thm}

\proof
We use Verma's expansion once again. In (\ref{Verma}) set
\begin{align}
&
 r=2,\quad s=1,\quad t=4 \quad u=4,\quad k=2,\label{general-rstuk}\\
 &
 w=1,\qquad x=q, \label{gen-1}\\
&
a_R=(qs/ab,qs/ad),\qquad b_S=q^2s^2/abcd,\qquad
c_T=(q^{-n},q^{n+1}s^2/a^2,s,qs/ac),\nonumber\\
&
d_U=(q^{1/2}s/a,-q^{1/2}s/a,,qs/a,-qs/a),\qquad e_K=(-qs/a^2,-qs/ac),\label{gen-2}\\
&
\gamma=qs^2/a^2c^2. \nonumber
\end{align}
The choice of the parameters (\ref{general-rstuk}) means that 
we have the expansion of the form ${}_6\phi_5=\sum {}_6\phi_5\cdot  {}_4\phi_3 $.
Note that the parameters in (\ref{gen-1}) and  (\ref{gen-2}) are
chosen in such a way that the structure of the resulting ${}_6\phi_5$ 
possesses the same structure as the one in 
$\Psi(x;a|a,-a,c,-c|q)$.

\begin{lem}
It holds that
\begin{align}
&
\Psi(x;s|a,b,c,d|q)\nonumber\\
&={(ax;q)_\infty\over (qx/a;q)_\infty}
\sum_{n\geq 0}
{(qs^2/a^2;q)_n \over (q;q)_n} (a x/s)^n \nonumber\\
&\times \sum_{j\geq 0}
{(q^{-n},q^{n+1}s^2/a^2,s,qs/ac,-qs/a^2,-qs/ac;q)_j \over
(q,q^{1/2}s/a,-q^{1/2}s/a,qs/a,-qs/a,q^{j+1} s^2/a^2c^2;q)_j} (-1)^j q^{j+\left(j\atop 2\right)} 
\label{gen-6}
\\
&\times
{}_6\phi_5
\left[
{q^{-n+j},q^{n+j+1}s^2/a^2,q^j s,q^{j+1}s/ac,-q^{j+1}s/a^2,-q^{j+1}s/ac 
\atop
q^{2j+2}s^2/a^2c^2,q^{1/2+j}s/a,-q^{1/2+j}s/a,q^{j+1}s/a,-q^{j+1}s/a
};q,q
\right]\nonumber\\
&\times 
{}_4\phi_3
\left[
{q^{-j},q^{j+1}s^2/a^2c^2,qs/ab,qs/ad
\atop
q^{2}s^2/abcd,-qs/a^2,-qs/ac
};q,q
\right].\nonumber
\end{align}
\end{lem}

Set $n=m+j$,
and change the order of the summation as $\sum_{n= 0}^\infty \sum_{j=0}^n=\sum_{j= 0}^\infty \sum_{m=0}^\infty$,
and apply Theorem \ref{thm:special} to the RHS of (\ref{gen-6}).
\begin{lem}
We have
\begin{align}
&\Psi(x;s|a,b,c,d|q)=
\sum_{j\geq 0} 
\sum_{k,l\geq 0}c_e(k,l;q^j s)x^{2k+2l}
{(qs^2/a^2;q)_{j}\over (q;q)_{j}} (ax/s)^{j}\label{general-5}\\
&\times
{(q^{-j},q^{j+1}s^2/a^2,s,qs/ac,-qs/a^2,-qs/ac;q)_j \over
(q,q^{1/2}s/a,-q^{1/2}s/a,qs/a,-qs/a,q^{j+1} s^2/a^2c^2;q)_j} (-1)^j q^{j+\left(j\atop 2\right)} \nonumber
\\
&\times 
{}_4\phi_3
\left[
{q^{-j},q^{j+1}s^2/a^2c^2,qs/ab,qs/ad
\atop
q^{2}s^2/abcd,-qs/a^2,-qs/ac
};q,q
\right]. \nonumber
\end{align}
\end{lem}

Applying the Sears' transformation and simplifying the factors, we have
\begin{align}
\mbox{RHS of (\ref{general-5})}&=
\sum_{j\geq 0} 
\sum_{k,l\geq 0}c_e(k,l;q^j s)x^{2k+2l} 
{(-d/c,qs/ab,s,qs^2/a^2c^2;q)_{j}\over 
(q,q^{2}s^2/abcd,q^{1/2}s/ac,-q^{1/2}s/ac;q)_{j}}(qx/d)^j\nonumber\\
&\times
{}_4\phi_3
\left[
{q^{-j},-b/a,qs/cd,-q^{-j}ac/s
\atop
-qs/ac,-q^{-j+1}c/d,q^{-j}ab/s
};q,q
\right]. \label{gen-7}
\end{align}

\begin{lem}
We have
\begin{align}
&\sum_{m,n\geq 0}c_o(m,n;s)x^{m+n}=\sum_{l\geq 0} x^{l}\sum_{m=0}^l c_o(m,l-m;s)\label{sum-co}\\
&=
\sum_{l\geq 0} 
{(-d/c,qs/ab,s,qs^2/a^2c^2;q)_{l}\over 
(q,q^{2}s^2/abcd,q^{1/2}s/ac,-q^{1/2}s/ac;q)_{l}}(qx/d)^l 
{}_4\phi_3\left[{q^{-l},-q^{-l}ac/s,-b/a,qs/cd \atop -q^{-l+1}c/d,q^{-l}ab/s,-qs/ac};q,q\right]. \nonumber
\end{align}
\end{lem}
\proof Streightforward calculation. \qed

Summarizing these, we have
\begin{align}
\mbox{RHS of (\ref{gen-7})}=
\sum_{m,n\geq 0}\sum_{k,l\geq 0} c_e(k,l;q^{m+n}s)c_o(m,n;s) x^{2k+2l+m+n}=\Phi(x;s|a,b,c,d|q),
\end{align}
thereby completing the proof of Theorem \ref{thm:general}. \qed

\section{Bibasic transformations  for $c_e(k,l;s)$ and $c_o(k,l;s)$}\label{BIBASIC}
In this section, we present some transformation formulas for
the series $\sum_{l=0}^k c_e(k-l,l;s)$ and $\sum_{l=0}^k c_o(k-l,l;s)$.
\subsection{Transformation by Watson's formula}

\begin{prp}
We have
\begin{align}
\sum_{l=0}^k c_e(k-l,l;s)=
{(a^2c^2/q,s^2;q^2)_k\over (q^2,q^3s^2/a^2c^2;q^2)_k}
(q^3/a^2c^2)^k
{}_4\phi_3\left[
{a^2/q,c^2/q,q^{2k}s^2,q^{-2k} \atop
-s,-qs,a^2c^2/q};q^2,q^2
\right]. \label{Watson}
\end{align}
\end{prp}
\begin{rmk}
This is manifestly symmetric in the exchange of  $a$ and $c$.
\end{rmk}

\proof
We have
\begin{align}
\sum_{l=0}^k {c_e(k-l,l;s)\over c_e(k,0;s)}
={}_8\phi_7
\left[
{s^2/a^2,q^2 s/a,-q^2 s/a,q^{2k} s^2,c^2/q,-q s/a^2,-q^2 s/a^2,q^{-2k} \atop
s/a,-s/a,q^{-2k+2}/a^2,q^3 s^2/a^2 c^2,-qs,-s,q^{2k+2}s^2/a^2}
;q^2,q^2/c^2
\right].
\end{align}
Then we apply Watson's transformation formula \cite[p.35,(2.5.1)]{GR} to the ${}_8\phi_7$ series.
Simplifying the factors, we have (\ref{Watson}).
\qed

\subsection{Bibasic transformation by Singh's formula}
\begin{prp}\label{prp:type-2}
We have
\begin{align}
\sum_{k,l\geq 0} c_e(k,l;s) x^{2k+2l}&=
\sum_{k,l\geq 0} 
{(qa^2/c^2,q^{2l}s^2;q^2)_k\over 
(q^2,q^{2l} q^3s^2/a^2c^2;q^2)_k} (q^2x^2/a^2)^k \label{type-2}\\
&\quad\quad\times 
{(c^2/q,s;q)_l (q^2s^2/a^4;q^2)_l\over 
(q,qs/a^2;q)_l (q^3 s^2/a^2c^2;q^2)_l} (q^2x^2/c^2)^l.\nonumber
\end{align}
\end{prp}

\proof
Applying Singh's formula and Sears' transformation to the ${}_4\phi_3$ series in  (\ref{Watson}), we have
\begin{align}
&{}_4\phi_3\left[
{a^2/q,c^2/q,q^{2k}s^2,q^{-2k} \atop
-s,-qs,a^2c^2/q};q^2,q^2
\right]=
{}_4\phi_3\left[
{a^2/q,c^2/q,q^{k}s,q^{-k} \atop
q^{-1/2}ac,-q^{-1/2}ac,-s};q,q
\right]\\
&={(qa^2/c^2;q^2)_k\over (a^2c^2/q;q^2)_k} (c^2/q)^k
{}_4\phi_3\left[
{q^{-k}-q^{-k},c^2/q,-qs/a^2 \atop
-s,q^{-k+1/2}c/a,-q^{-k+1/2}c/a};q,q
\right]. \nonumber
\end{align}
Simplifying the factors, we have (\ref{type-2}).
\qed

\subsection{Bibasic transformation by $q$-analogue of Bailey's formula}
\begin{prp}\label{prp:type-3}
We have
\begin{align}
\sum_{k,l\geq 0} c_e(k,l;s) x^{2k+2l}&=
\sum_{k,l\geq 0} 
{(q a^2/c^2;q^2)_k (s;q^2)_k(q s;q^2)_k (q^2 s^2/c^4;q^2)_k \over 
(q^2;q^2)_k(q s/c^2;q^2)_k(q^2 s/c^2;q^2)_k(q^3 s^2/a^2 c^2;q^2)_k}(q^2x^2/a^2)^k \nonumber \\
&\times
{(c^2/q;q)_l (q^{2k} s;q)_l \over (q;q)_l (q^{2k} q^2 s/c^2;q)_l }(q^2x^2/c^2)^l\nonumber \\
&=
\sum_{k,l\geq 0} 
{(q a^2/c^2;q^2)_k (q^3 s/c^2;q^2)_k (q^2 s^2/c^4;q^2)_k \over 
(q^2;q^2)_k(q s/c^2;q^2)_k(q^3 s^2/a^2 c^2;q^2)_k}(q^2x^2/a^2)^k \nonumber \\
&\times
{(c^2/q;q)_l (s;q)_{2k+l} \over (q;q)_l (q^2 s/c^2;q)_{2k+l} }(q^2x^2/c^2)^l. \label{type-3}
\end{align}
\end{prp}

\proof
The coefficient of $x^{2l}$ on RHS of (\ref{type-3}) reads
\begin{align}
&
{(c^2/q;q)_l ( s;q)_l \over (q;q)_l (q^2 s/c^2;q)_l }(q^2/c^2)^l\\
&\times \sum_{k=0}^l
{(q^2s^2/c^4,q^3s/c^2,qa^2/c^2;q^2)_k\over 
(q^2,qs/c^2,q^3s^2/a^2c^2;q^2)_k}
{(q^ls,q^{-l};q)_k\over 
(q^{-l+2}/c^2,q^{l+2}s/c^2;q)_k}(q^2/a^2)^k. \nonumber
\end{align}
The second factor can be expressed in terms of the bibasic
hypergeometric series as
\begin{align}
\Phi\left[
{q^2s^2/c^4,q^3s/c^2,qa^2/c^2:q^ls,q^{-l} \atop
qs/c^2,q^3s^2/a^2c^2: q^{-l+2}/c^2,q^{l+2}s/c^2};q^2,q\,;q^2/a^2
\right], \label{Phi-bibasic}
\end{align}
where we used the notation for the bibasic
hypergeometric series (see \cite[p.85, (3.9.1)]{GR})
\begin{align}
\Phi\left[ 
{a_1,\ldots,a_{r+1}:c_1,\ldots,c_s \atop
b_1,\ldots,b_r:d_1,\ldots,d_s};q,p;z
\right]
=\sum_{n=0}^\infty
{(a_1,\ldots,a_{r+1};q)_n\over (q,b_1,\ldots,b_r;q)_n}
{(c_1,\ldots,c_s;p)_n\over (d_1,\ldots,d_s;p)_n} z^n.
\end{align}

Recall the $q$-analogue of Bailey's transformation \cite[p.89, (3.10.14)]{GR}
\begin{align}
&
\Phi\left[
{a^2,aq^2,-aq^2,b^2,c^2:-aq/w,q^{-n}\atop
a,-a,a^2q^2/b^2,a^2q^2/c^2:w,-a q^{n+1}};q^2,q\,;{awq^{n+1}\over b^2 c^2}\right] \label{Baily}\\
&=
{(w/a,-aq;q)_n\over (w,-q;q)_n}
{(wq^{-n-1}/a,a q^{2-n}/w;q^2)_n\over 
(a q^{1-n}/w,w q^{-n}/a;q^2)_n} 
{}_5\phi_4\left[
{aq,aq^2,a^2q^2/b^2c^2,a^2q^2/w^2,q^{-2n}
\atop
a^2q^2/b^2,a^2q^2/c^2,aq^{2-n}/w,aq^{3-n}/w};q^2,q^2
\right]. \nonumber
\end{align}
Note that setting $a=c^2$ 
and replacing the parameters as 
$(a,w,b^2)\rightarrow (-qs/c^2,q^{-l+2}/c^2,qa^2/c^2)$
in (\ref{Baily}) we have a bibasic
transformation formula involving a ${}_4\phi_3$ series which applies to (\ref{Phi-bibasic}).
The resulting ${}_4\phi_3$ has the structure
\begin{align}
{}_4\phi_3\left[
{-q^2s/a^2,-q^2s/c^2,q^{2l}s^2,q^{-2l} \atop
-s,-qs,q^3s^2/a^2c^2};q^2,q^2
\right],
\end{align}
which is transformed by Sears' formula into 
\begin{align}
{}_4\phi_3\left[
{a^2/q,c^2/q,q^{2l}s^2,q^{-2l} \atop
-s,-qs,a^2c^2/q};q^2,q^2
\right],
\end{align}
up to a multiplicative factor given in terms of $q$-shifted factorials. Summarizing these,
one finds that RHS of (\ref{type-3}) is transformed to RHS of (\ref{Watson}). 

\subsection{Simplification of $\Phi(x;s|-a,b,-q^{1/2} a,q^{1/2} b|q)$}
One finds a simplification of the $\Phi(x;s)$ 
to a twofold series with bases $q^{1/2}$ and $q$
when we specialize the parameters as
$(a,b,c,d)\rightarrow (-a,b,-q^{1/2} a,q^{1/2} b)$.

\begin{prp}\label{simplify-1}
We have
\begin{align}
&\Phi(x;s|-a,b,-q^{1/2} a,q^{1/2} b|q)=
\sum_{l,m\geq 0}
{(a^2;q)_m(q^l s;q)_m\over (q;q)_m(q^{l+1}s/a^2;q)_m} (q x^2/a^2)^m\\
&\quad \quad \times
{(b/a;q^{1/2})_l(s/a^2;q^{1/2})_l  (s;q)_l\over
 (q^{1/2};q^{1/2})_l (q^{1/2}s/ab;q^{1/2})_l (s/a^2;q)_l }(q^{1/2}x/b)^l. \nonumber
\end{align}
\end{prp}

\proof The ${}_4\phi_3$ series on the RHS of (\ref{sum-co}) 
written for the parameters $(-a,b,-c,d)$
is transformed by Sears' formula as
\begin{align}
{}_4\phi_3\left[
{q^{-l},-q^{-l}ac/s,b/a,-qs/cd \atop
q^{-l+1}c/d,-q^{-l}ab/s,-qs/ac};q,q
\right]
&=
{(bd/ac,-qs/bc;q)_l\over (d/c,-qs/ac;q)_l}
{}_4\phi_3\left[
{q^{-l},q^{-l-1}abcd/s^2,b/a,b/c \atop
-q^{-l}ab/s,-q^{-l}bc/s,bd/ac};q,q
\right]. \nonumber
\end{align}
Note that if $c=q^{1/2}a,d=q^{1/2}b$, Singh's transformation applies and we 
can simplify the ${}_4\phi_3$ series as
\begin{align}
&{}_4\phi_3\left[
{q^{-l},q^{-l}a^2b^2/s^2,b/a,q^{-1/2}b/a \atop
-q^{-l}ab/s,-q^{-l+1/2}ab/s,b^2/a^2};q,q
\right]
={}_4\phi_3\left[
{q^{-l/2},q^{-l/2}ab/s,b/a,q^{-1/2}b/a \atop
b/a,-b/a,-q^{-l}ab/s};q^{1/2},q^{1/2}
\right]\\
&=
{}_3\phi_2\left[
{q^{-l/2},q^{-l/2}ab/s,q^{-1/2}b/a \atop
-b/a,-q^{-l}ab/s};q^{1/2},q^{1/2}
\right]
=
{(-q^{1/2},q^{l/2}s/a^2;q^{1/2})_l\over (-b/a,-q^{l/2+1/2}s/ab;q^{1/2})_l},\nonumber
\end{align}
where we used the $q$-Saalsch\"utz sum \cite[p.813, (1.7.2)]{GR}.
Hence we have
\begin{align}
\sum_{m=0}^l c_o(m,l-m;s|-a,b,-q^{1/2} a,q^{1/2} b,q)=
{(b/a;q^{1/2})_l(s/a^2;q^{1/2})_l  (s;q)_l\over
 (q^{1/2};q^{1/2})_l (q^{1/2}s/ab;q^{1/2})_l (s/a^2;q)_l }(q^{1/2}/b)^l.
 \end{align}
{} From (\ref{type-2}) or (\ref{type-3}), we have
\begin{align}
\sum_{k=0}^m c_e(k,m-k;q^l s|-a,b,-q^{1/2} a,q^{1/2} b,q)=
{(a^2;q)_m(q^l s;q)_m\over (q;q)_m(q^{l+1}s/a^2;q)_m} (q /a^2)^m.
 \end{align}
\qed

\subsection{Simplification of $\Phi(x;s|-a,b,-q^{1/2} a,q^{1/2} a|q)$}
One finds another simplification of the $\Phi(x;s)$ 
to a twofold series with base $q$
when we specialize the parameters as
$(a,b,c,d)\rightarrow (-a,b,-q^{1/2} a,q^{1/2} a)$.

\begin{prp}\label{simplify-2}
We have
\begin{align}
&\Phi(x;s|-a,b,-q^{1/2} a,q^{1/2} a|q)=
\sum_{l,m\geq 0}
{(a^2;q)_m(q^l s;q)_m\over (q;q)_m(q^{l+1}s/a^2;q)_m} (q x^2/a^2)^m\\
&\quad \quad \times
{(b/a;q)_l(s^2/a^4;q)_l  (s;q)_l\over
 (q;q)_l (qs^2/a^3b;q)_l (s/a^2;q)_l }(qx/b)^l. \nonumber
\end{align}
\end{prp}
\proof
Simple calculation using (\ref{co}) and (\ref{type-2}) (or (\ref{type-3})).

\section{Proof of Theorem \ref{thm:g_r}}\label{Koornwinder}
We apply our formulas for $f(x;s)$ for Koornwinder polynomials with one row diagram.
We use the kernel function of type $BC$ which intertwines the action of the
Koornwinder operator of type $BC_n$ and type $BC_1$.
As for the notations and basic facts about the kernel function of type $BC$, 
we refer the readers to \cite{KNS} and to Section \ref{appendix}. 

\begin{lem}
{} Note that from Proposition \ref{prp:type-3},
we have
\begin{align}
&(1-x^2)\sum_{k,l\geq 0} c_e(k,l;s) x^{2k+2l}=
\sum_{k,l\geq 0}  
c'_e(k,l;s|a,c|q) x^{2k+2l}, 
\end{align}
where
\begin{align}
&
c'_e(k,l;s|a,c|q)= 
{(q a^2/c^2;q^2)_k (q^3 s/c^2;q^2)_k (q^2 s^2/c^4;q^2)_k \over 
(q^2;q^2)_k(q s/c^2;q^2)_k(q^3 s^2/a^2 c^2;q^2)_k}(q^2x^2/a^2)^k \nonumber \\
&\quad \quad \times
{(c^2/q^2;q)_l (s/q;q)_{2k+l} \over (q;q)_l (q^2 s/c^2;q)_{2k+l} }
{1-q^{2k+2l-1}s\over 1-q^{-1}s}
(q^2x^2/c^2)^l. \label{type-4} 
\end{align}
\end{lem} 
\proof Streightforward. \qed

\subsection{Koornwinder polynomial with one row diagram $P_{(r)}(x|a,b,c,d|q,t)$} \label{P(r)}
We move on to the proof of Theorem \ref{thm:g_r}.
Recall  that $n$  is a positive integer, $x=(x_1,\ldots,x_n)$ is a set of variables, and 
$P_{(r)}(x|a,b,c,d|q,t)$ denotes the Koornwinder polynomial with one row diagram $(r)$.

\noindent
{\it Proof of Theorem \ref{thm:g_r}}.
We consider the following special case of Theorem \ref{thm:reproduction} below
\begin{align}
&x=(x_1,\ldots,x_n) \quad (n\in \bbZ_{>0}),\\
&y=(y_1) \quad \quad\quad\quad\,\, (m=1),\\
&
\Pi(x;y)=y^{\beta n}\prod_{i=1}^n
{(q^{1/2}t^{1/2}yx_i;q)_\infty \over (q^{1/2}t^{-1/2}yx_i;q)_\infty }
{(q^{1/2}t^{1/2}y/x_i;q)_\infty \over (q^{1/2}t^{-1/2}y/x_i;q)_\infty },\\
&\lambda=(r)\quad \quad\quad\quad\,\, \,\, (r\in \bbZ_{\geq 0}),\\
&s=(s_1)=q^{1-r}t^{-n},\\
&V(y)=y^{-1}(1-y^2),\\
&\widehat{f}(x;s)=y^{-r+1-n\beta} \widehat{\Phi}(x;s)=
y^{-r+1-n\beta} \sum_{k,l,i,j\geq 0}
\widehat{c}_e(k,l;q^{i+j} s)\widehat{c}_o(i,j; s) x^{2k+2l+i+j},
\end{align}
where 
\begin{align}
&\widehat{c}_e(k,l;s)=c_e(k,l;s|\sqrt{q/t}a,\sqrt{q/t}c|q),\\
&\widehat{c}_o(i,j; s)=c_o(i,j;s|\sqrt{q/t}a,\sqrt{q/t}b,\sqrt{q/t}c,\sqrt{q/t}d|q).
\end{align}
Note that we have
\begin{align}
y^{-n\beta}\Pi(x;y)=\sum_{r\geq 0}G_r(x;q,t) (q/t)^{r/2} y^r,
\end{align}
where $G_r$ is defined in Definition \ref{G_r}.

Then we calculate the constant term, which is proportional to $P_{(r)}(x|q,b,c,d|q,t)$,  as
\begin{align}
&
\biggl[  \Pi(x;y)V(y) \widehat{f}(x;s)\biggr]_{1,y}\\
&=
\biggl[\sum_{\theta\geq 0}G_\theta (x;q,t) (q/t)^{\theta /2} y^\theta  y^{-r} (1-y^2) \widehat{\Phi}(x;s)\biggr]_{1,y}\nonumber \\
&=
\sum_{k,l,i,j\geq 0} 
G_{r-2k-2l-i-j} (x;q,t) (q/t)^{r-2k-2l-i-j}
\widehat{c}'_e(k,l;q^{1-r+i+j} t^{-n})\widehat{c}_o(i,j; q^{1-r}t^{-n}) \nonumber \\
&= (q/t)^{r/2 }{(t;q)_r\over (q;q)_r} P_{(r)}(x|q,b,c,d|q,t).
\end{align}
where $\widehat{c}'_e(k,l;s)= c'_e(k,l;s|\sqrt{q/t}a,\sqrt{q/t}c|q)$ (see (\ref{type-4}) above).  
Here we have used 
\begin{align}
G_r(x;q,t)={(t;q)_r\over (q;q)_r}  \sum_{i=1}^n x_i^n + \mbox{lower degree terms}.
\end{align}
\qed

\subsection{Macdonald polynomial of type $C_n$ with one row diagram $P^{(C_n)}_{(r)}(x|b;q,t)$}\label{Mac-C}
Let $P^{(C_n)}_{\lambda}(x|b;q,t)$ be the Macdonald polynomial of type $C_n$.
(As for the notation, see \S\S \ref{notation-C}). In view of (\ref{param-C}), 
we need $\Phi(x;s|a,b,c,d|q)$ written for the parameters
\begin{eqnarray}
(-(qb/t)^{1/2}, (qb/t)^{1/2},-q^{1/2}(qb/t)^{1/2},q^{1/2}(qb/t)^{1/2},q),
\end{eqnarray}
for the construction of the formula for $P^{(C_n)}_{(r)}(x|q,t)$.

\begin{lem}\label{type-C}
{}From (\ref{type-2}) or (\ref{type-3}), we have
\begin{align}
&(1-x^2) \Phi(x;s|-a,a,-q^{1/2}a,q^{1/2}a|q)=
(1-x^2)\sum_{j\geq 0}
{(a^2;q)_j (s;q)_j\over (q;q)_j (qs/a^2;q)_j} (qx^2/a^2)^j\nonumber \\
&\quad =
\sum_{j\geq 0}
{(a^2/q;q)_j (s/q;q)_j\over (q;q)_j (qs/a^2;q)_j}{1-q^{2j-1}s\over 1-q^{-1}s} (qx^2/a^2)^j.
\end{align}
\end{lem}

\begin{thm}\label{Thm-C}
We have
\begin{align}
 P^{(C_n)}_{(r)}(x)=
 {(q;q)_r\over (t;q)_r}
 \sum_{j=0}^{\lfloor r/2\rfloor}G_{r-2j}(x;q,t)
{(b/t;q)_j(t^{-n}q^{-r};q)_j\over
(q;q)_j(t^{-n+1}q^{-r+1}/b;q)_j}{1-t^{-n}q^{-r+2j} \over
 1-t^{-n}q^{-r}}(t^2/qb)^j,\nonumber
\end{align}
thereby proving Lassalle's conjecture for type $C$ \cite[p.8, Conjecture 1]{Lass} (See \S\S\ref{Lassalle-con}).
\end{thm}

\subsection{Macdonald polynomial of type $B_n$, $D_n$ 
with one row diagram $P^{(B_n)}_{(r)}(x|a;q,t)$, $P^{(D_n)}_{(r)}(x|q,t)$}\label{Mac-BD}

For the construction of the formula for $P^{(B_n)}_{(r)}(x|a;q,t)$, 
we need $\Phi(x;s|a,b,c,d|q)$ written for the parameters
\begin{eqnarray}
(-(q/t)^{1/2}, a(q/t)^{1/2},-q t^{-1/2},qt^{-1/2};q,q/t).
\end{eqnarray}

\begin{lem}\label{type-B}
We have
\begin{align}
&(1-x^2)\Phi(x;s|-a, b,-q^{1/2}a ,q^{1/2}a|q)\\
&=
(1-x^2)\sum_{i,j\geq 0}
{(a^2;q)_j(s;q)_{i+j}\over (q;q)_j(qs/a^2;q)_{i+j}} 
{(b/a;q)_i(s^2/a^4;q)_i(qs/a^2;q)_i \over
 (q;q)_i (qs^2/a^3 b;q)_i (s/a^2;q)_i}(q x^2/a^2)^j(qx/b)^i \nonumber\\
 &=
 \sum_{i,j\geq 0}
{(a^2/q;q)_j(s/q;q)_{i+j}\over (q;q)_j(qs/a^2;q)_{i+j}} 
{(b/a;q)_i(s^2/a^4;q)_i(qs/a^2;q)_i \over
 (q;q)_i (qs^2/a^3 b;q)_i  (s/a^2;q)_i}
 {1-q^{i+2j-1}s\over 1-q^{-1}s}
 (q x^2/a^2)^j(qx/b)^i. \nonumber
\end{align}
\end{lem}


\begin{thm}\label{Thm-BD}
We have
\begin{align}
&P^{(B_n)}_{(r)}(x|a;q,t)=
{(q;q)_r\over (t;q)_r}
\sum_{0\leq i+2j\leq r}G_{r-i-2j}(x;q,t)\\
&\times {(a;q)_{i} 
( t^{-n+1}q^{-r+1};q)_{i}
( t^{-n}q^{-r};q)_{i+j}
( t^{-2n+2}q^{-2r};q)_{i} (1/t;q)_{j}
\over
 (q;q)_{i}
 (t^{-n+1}q^{-r};q)_{i}
 (t^{-n+1}q^{-r+1};q)_{i+j} 
 (t^{-2n+2}q^{-2r+1}/a;q)_{i}  (q;q)_{j}}
{1-t^{-n}q^{-r+i+2j} \over 
1- t^{-n}q^{-r}}
(t/a)^{i} (t^2/q)^{j}. \nonumber
\end{align}
By setting $a=1$ we have
\begin{align}
&P^{(D_n)}_{(r)}(x|q,t)=
{(q;q)_r\over (t;q)_r}
\sum_{0\leq 2j\leq r}G_{r-2j}(x;q,t)
{
( t^{-n}q^{-r};q)_{j}
(1/t;q)_{j}
\over
 (t^{-n+1}q^{-r+1};q)_{j} 
   (q;q)_{j}}
{1-t^{-n}q^{-r+2j} \over 
1- t^{-n}q^{-r}}
(t^2/q)^{j}. \nonumber
\end{align}
Hence we have proved Lassalle's conjecture for type $D$ and $B$  \cite[p.10, Conjecture 3 and p.11, Conjecture 4 ]{Lass} 
(See \S\S\ref{Lassalle-con}).
\end{thm}

\section{Appendix}\label{appendix}
We recall briefly some properties concerning the Koornwinder polynomials
needed for the construction given in Section \ref{Koornwinder}.

\subsection{Kernel function $\Pi(x;y)$}\label{Kernel}

Let $(a,b,c,d;q,t)$ be a set of complex parameters with $|q|<1$. 
Set  $\alpha=(abcdq^{-1})^{1/2}$ for simplicity. 
Let $x=(x_1,\ldots,x_n)$ be a set of independent indeterminates. 
Koornwinder's $q$-difference operator  ${\mathcal D}_x={\mathcal D}_x(a,b,c,d|q,t)$ is defined by \cite{Ko}
\begin{eqnarray}
&&{\mathcal D}_x=
\sum_{i=1}^n {(1-ax_i)(1-bx_i)(1-cx_i)(1-dx_i)\over 
\alpha t^{n-1}(1-x_i^2)(1-qx_i^2)}
\prod_{j\neq i} {(1-t x_ix_j)(1-t x_i/x_j)\over (1-x_ix_j)(1-x_i/x_j)}
\left(T_{q,x_i}-1\right) \\
&&+
\sum_{i=1}^n {(1-a/x_i)(1-b/x_i)(1-c/x_i)(1-d/x_i)\over 
\alpha t^{n-1}(1-1/x_i^2)(1-q/x_i^2)}
\prod_{j\neq i} {(1-t x_j/x_i)(1-t /x_ix_j)\over (1-x_j/x_i)(1-1/x_ix_j)}
\left(T_{q^{-1},x_i}-1\right),\nonumber
\end{eqnarray}
where  $T_{q,x}^{\pm1}f(x_1,\ldots,x_i,\ldots ,x_n)=f(x_1,\ldots,q^{\pm 1}x_i,\ldots ,x_n)$. 
Koornwinder polynomial $P_\lambda(x)=P_\lambda(x|a,b,c,d|q,t)$ 
with partition $\lambda=(\lambda_1,\ldots,\lambda_n)$
(i.e. $\lambda_i\in \bbZ_{\geq 0},\lambda_1\geq \cdots\geq \lambda_n$)
is uniquely characterized by the two conditions
(a) $P_\lambda(x)$ is a $S_n \ltimes (\bbZ/2\bbZ)^n$ invariant Laurent
polynomial having the triangular expansion in terms of the monomial 
basis $(m_\lambda)$ as $P_\lambda(x)=m_\lambda(x)+\mbox{lower terms}$,
(b) $P_\lambda(x)$ satisfies ${\mathcal D}_x P_\lambda=d_\lambda P_\lambda$.
The eigenvalue is given by
\begin{align}
&&d_\lambda=\sum_{j=1}^n \langle abcdq^{-1}t^{2n-2j}q^{\lambda_j}\rangle
\langle q^{\lambda_j}\rangle
=
\sum_{j=1}^n \langle \alpha t^{n-j}q^{\lambda_j}; \alpha t^{n-j}
\rangle,
\end{align}
where we used the notations 
$\langle x\rangle=x^{1/2}-x^{-1/2}$ and
$\langle x;y\rangle=\langle xy\rangle\langle x/y\rangle=x+x^{-1}-y-y^{-1}$
for simplicity of display.

\begin{dfn}
Define the involution $\widetilde{*}$ of the parameters by
\begin{align}
\widetilde{a}={\sqrt{qt}\over a},\qquad 
\widetilde{b}={\sqrt{qt}\over b},\qquad 
\widetilde{c}={\sqrt{qt}\over c},\qquad 
\widetilde{d}={\sqrt{qt}\over d},\qquad 
\widetilde{q}=q,\qquad 
\widetilde{t}=t.\end{align}
We write $\widetilde{\mathcal D}_x={\mathcal D}_x(
\widetilde{a},\widetilde{b},\widetilde{c},\widetilde{d}|\widetilde{q},\widetilde{t})$
and $\widetilde{P}_\lambda(x)=P_\lambda(x|
\widetilde{a},\widetilde{b},\widetilde{c},\widetilde{d}|\widetilde{q},\widetilde{t}))$ for simplicity of display.
\end{dfn}

\begin{prp}
Let  $n$ and $m$ be positive integers, and let
$x=(x_1,\cdots,x_n)$, $y=(y_1,\cdots,y_m)$ be two sets of independent indeterminates.
Let $\beta\in \bbC$ be satisfying $t=q^\beta$.
Set
\begin{align}
\Pi(x;y)=\prod_{k=1}^m y_k^{\beta n}\prod_{i=1}^n \prod_{j=1}^m
{(q^{1/2}t^{1/2}y_jx_i;q)_\infty \over (q^{1/2}t^{-1/2}y_jx_i;q)_\infty }
{(q^{1/2}t^{1/2}y_j/x_i;q)_\infty \over (q^{1/2}t^{-1/2}y_j/x_i;q)_\infty }. \label{Kernel-function}
\end{align}
Then we have the kernel function identity \cite{KNS}
\begin{align}
{\mathcal D}_x \Pi(x;y)-\widetilde{{\mathcal D}}_y \Pi(x;y)=
-{1\over \langle t\rangle} \langle t^n\rangle\langle t^{m}\rangle
\langle abcdq^{-1}t^{n-m-1}\rangle \Pi(x;y).
\end{align}
\end{prp}

\subsection{Series $f(x;s)$ for $BC_n$}
Let $s=(s_1,\ldots,s_n)$ be a set of complex parameters. 
Corresponding to $s$, 
we introduce $\lambda=(\lambda_1,\ldots,\lambda_n)$
by the conditions $s_i=t^{-n+i}q^{-\lambda_i}$ ($i=1,\ldots,n$). 
We use the notation for the multiple index as $x^{-\lambda}=\prod_ix_i^{-\lambda_i}$. 
Let $f(x;s)\in x^{-\lambda} \bbC[[x_1/x_2,\ldots,x_{n-1}/x_n,x_n]]$ be the infinite series 
satisfying the conditions
\begin{align}
&f(x;s)=x^{-\lambda}\sum_{\alpha \in Q^+} c_\alpha(s)x^\alpha,\quad c_0(s)=1,\\
&{\mathcal D}_x f(x;s)= \sum_{i=1}^n \langle \alpha s_i^{-1} ;\alpha t^{n-1} \rangle  f(x;s).
\end{align}

\subsection{Reproduction formula}\label{reproduction}
Let ${\mathcal D}_x^*$ be the adjoint of ${\mathcal D}_x$ given by  
\begin{align}
&{\mathcal D}_x^*=
\sum_{i=1}^n 
\left(T_{q,x_i}^{-1}-1\right)
{(1-ax_i)(1-bx_i)(1-cx_i)(1-dx_i)\over 
\alpha t^{n-1}(1-x_i^2)(1-qx_i^2)}
\prod_{j\neq i} {(1-t x_ix_j)(1-t x_i/x_j)\over (1-x_ix_j)(1-x_i/x_j)}
 \\
&+
\sum_{i=1}^n 
\left(T_{q,x_i}-1\right){(1-a/x_i)(1-b/x_i)(1-c/x_i)(1-d/x_i)\over 
\alpha t^{n-1}(1-1/x_i^2)(1-q/x_i^2)}
\prod_{j\neq i} {(1-t x_j/x_i)(1-t /x_ix_j)\over (1-x_j/x_i)(1-1/x_ix_j)}. \nonumber
\end{align}
Denote by $V(x)$ the Weyl denominator of type 
$BC_n$
\begin{align}
V(x)=
\prod_{k=1}^n x_k^{-n+k-1} \prod_{i=1}^n (1-x_i^{2})
\prod_{1\leq i<j\leq n}(1-x_ix_j)(1-x_i/x_j).
\end{align}

\begin{dfn}
Define the involution $\overline{*}$ of the parameters by
\begin{align}
\overline{a}=q/a,\qquad 
\overline{b}=q/b,\qquad 
\overline{c}=q/c,\qquad 
\overline{d}=q/d,\qquad 
\overline{q}=q,\qquad
\overline{t}=q/t.
\end{align}
Write for simplicity the composition of the involutions as $\widehat{*}=\overline{\widetilde{*}}$, namely we have
\begin{align}
\widehat{a}=\sqrt{q/t}a,\qquad 
\widehat{b}=\sqrt{q/t}b,\qquad 
\widehat{c}=\sqrt{q/t}c,\qquad 
\widehat{d}=\sqrt{q/t}d,\qquad 
\widehat{q}=q,\qquad
\widehat{t}=q/t.
\end{align}
\end{dfn}

\begin{prp}
We have
\begin{align}
V(x)^{-1}{\mathcal D}_x^*V(x)-
 \overline{\mathcal D}_x
 =\sum_{j=1}^n \langle\overline{\alpha}\bar{t}^{n-j};\alpha t^{n-j}\rangle.
 \end{align}
\end{prp}

\begin{thm}\label{thm:reproduction}
Let $n\geq m$ be positive integers, and 
$x=(x_1,\ldots,x_n), y=(y_1,\ldots,y_m)$ be sets of variables.
Let $\lambda=(\lambda_1,\ldots,\lambda_m)$ be a partition of length $\leq m$.
Set
\begin{align}
s_i=\widehat{t}^{-m+i}q^{-\lambda_{m+1-i} +m+1-i-n \beta}\quad (1\leq i\leq m).\label{choice-s}
\end{align}
Let $\widehat{f}(y;s)$ denotes the formal series in $y$ characterized by
\begin{align}
&\widehat{f}(y;s)=\prod_{i=1}^m y_i^{-\lambda_{m+1-i} +m+1-i-n \beta}\sum_{\alpha\in Q^+}
\widehat{c}_\alpha(s) y^\alpha,\\
&
\widehat{\mathcal D}_y \widehat{f}(y;s)=\sum_{i=1}^m \langle \widehat{\alpha }s_i^{-1};\widehat{\alpha} 
\widehat{t}^{m-i}\rangle  \widehat{f}(y;s).
\end{align}
Then $\Pi(x;y)V(y) \widehat{f}(y;s)$ has no fractional powers in $y$, allowing us to consider the constant term in $y$. 
We have
\begin{align}
{\mathcal D}_x \biggl[ \Pi(x;y)V(y) \widehat{f}(y;s)\biggr]_{1,y}=
\sum_{i=1}^m \langle \alpha t^{n-1}q^{\lambda_i};\alpha t^{n-i}\rangle \biggl[ \Pi(x;y)V(y) \widehat{f}(y;s)\biggr]_{1,y},
\end{align}
where the symbol $[\cdots]_{1,y}$ denotes the constant term in $y$.
Hence $[\Pi(x;y)V(y) \widehat{f}(y;s)]_{1,y}$ gives us the Koormwinder polynomila $P_\lambda(x)$
up to a multiplication constant.
\end{thm}

\proof
Note that from the choice of $s$ (\ref{choice-s}) and the definition of $\Pi(x;y)$ in (\ref{Kernel-function}),
fractional powers in $y$ cancels in the combination $\Pi(x;y)V(y) \widehat{f}(y;s)$. We have
\begin{align}
&\biggl(
{\mathcal D}_x +
{\langle t^n\rangle \langle t^m\rangle \langle abcdq^{-1}t^{n-m-1}\rangle \over \langle t\rangle}
\biggr)
\biggl[ \Pi(x;y)V(y) \widehat{f}(y;s)\biggr]_{1,y}\nonumber \\
&= \biggl[  \biggl(\widetilde{{\mathcal D}}_y\Pi(x;y)\biggr)V(y) \widehat{f}(y;s)\biggr]_{1,y}\\
&= \biggl[ \Pi(x;y) \biggl(\widetilde{{\mathcal D}}_x ^* V(y) \widehat{f}(y;s) \biggr)\biggr]_{1,y}\nonumber \\
&= \biggl[ \Pi(x;y)  V(y)\biggl(
\biggl(
\widehat{{\mathcal D}}_x 
+\sum_{i=1}^m \langle \widehat{\alpha} \widehat{t}^{m-i};\widetilde{\alpha} t^{m-i} 
\rangle\biggr)\widehat{f}(y;s) \biggr)\biggr]_{1,y}.\nonumber
\end{align}
To calculate the eigenvalue, we need a
\begin{lem}
Write $\alpha=(abcdq^{-1})^{1/2}$, $\widetilde{\alpha}=t/\alpha$. We have
\begin{eqnarray}
{1\over \langle t\rangle} \langle t^n\rangle\langle t^{m}\rangle
\langle abcdq^{-1}t^{n-m-1}\rangle
=\sum_{i=1}^{m\wedge n} 
\langle \alpha t^{n-i};\widetilde{\alpha}t^{m-i}\rangle.
\end{eqnarray}
\end{lem}
Hence by noting $\widehat{\alpha}= q\alpha/t $ and (\ref{choice-s}) we have
\begin{align}
&\sum_{i=1}^m \langle \widehat{\alpha }s_i^{-1};\widehat{\alpha} \widehat{t}^{m-i}\rangle +
\sum_{i=1}^m \langle \widehat{\alpha} \widehat{t}^{m-i};\widetilde{\alpha} t^{m-i} \rangle-
{1\over \langle t\rangle} \langle t^n\rangle\langle t^{m}\rangle
\langle abcdq^{-1}t^{n-m-1}\rangle\\
=&
\sum_{i=1}^m \langle \widehat{\alpha }s_i^{-1}; \alpha t^{n-i}\rangle
=\sum_{i=1}^m \langle \alpha t^{n-m-1+i}q^{\lambda_{m+1-i}}; \alpha t^{n-i}\rangle
=\sum_{i=1}^m \langle \alpha t^{n-i}q^{\lambda_{i}}; \alpha t^{n-i}\rangle. \nonumber
\end{align}
\qed

\subsection{Macdonald polynomials of type $C$}\label{notation-C}
We consider some degeneration of the Koornwinder polynomials to
Macdonald polynomials. As for the details, we refer the readers to \cite{Ko} and \cite{Mac1}.
Setting the parameters as 
$(a,b,c,d,q,t)\rightarrow (-b^{1/2},ab^{1/2}, -q^{1/2}b^{1/2},q^{1/2}ab^{1/2},q,t)$ 
in the Koornwinder polynomial $P_\lambda(x)$, 
we obtain the Macdonald polynomials of type $(BC_n,C_n)$ \cite{Ko}
\begin{align}
P^{(BC_n,C_n)}_\lambda(x|a,b;q,t)=P_\lambda(x|-b^{1/2},ab^{1/2}, -q^{1/2}b^{1/2},q^{1/2}ab^{1/2}|q,t)
\end{align} 
Namely, setting
\begin{align}
D^{(BC_n,C_n)}_x=\sum_{\sigma_1,\cdots,\sigma_n=\pm 1}
\prod_{i=1}^n
{(1-a b^{1/2} x_i^{\sigma_i})(1+b^{1/2} x_i^{\sigma_i})\over 1-x_i^{2\sigma_i}}
\prod_{1\leq i<j\leq n}
{1-t x_i^{\sigma_i} x_i^{\sigma_j}\over 1-x_i^{\sigma_i} x_i^{\sigma_j}}
T_{q^{\sigma_1/2},x_1}\cdots T_{q^{\sigma_n/2},x_n},
\end{align}
we have 
\begin{align}
&P^{(BC_n,C_n)}_\lambda(x)=m_\lambda+\mbox{lower terms},\\
&D^{(BC_n,C_n)}_x P^{(BC_n,C_n)}_\lambda(x)=
(ab)^{n/2}t^{n(n-1)/4}\sum_{\sigma_1,\cdots,\sigma_n=\pm 1} 
 s_1^{\sigma_1/2}\cdots s_n^{\sigma_n/2}
\cdot 
P^{(BC_n,C_n)}_\lambda(x),
\end{align}
where $s_i=ab t^{n-i}q^{\lambda_i}$.

The special case $a=1$ is called the Macdonald polynomials of type $C_n$
\begin{align}
P^{(C_n)}_\lambda(x|b;q,t)=P^{(BC_n,C_n)}_\lambda(x|1,b;q,t).
\end{align}

Note that the application of the twist $\widehat{*}$ on the parameters 
$ (-b^{1/2},ab^{1/2}, -q^{1/2}b^{1/2},q^{1/2}ab^{1/2},q,t)$ gives
\begin{eqnarray}
(-(qb/t)^{1/2},a (qb/t)^{1/2},-q^{1/2}(qb/t)^{1/2},q^{1/2}a(qb/t)^{1/2};q,q/t). \label{param-C}
\end{eqnarray}

\subsection{Macdonald polynomials of type $B$ and $D$}\label{notation-BD}
Setting the parameters as 
$(a,b,c,d,q,t)\rightarrow (-b^{1/2},ab^{1/2},-q^{1/2},q^{1/2};q,t)$ 
in the Koornwinder polynomial $P_\lambda(x)$, 
we obtain the Macdonald polynomials of type $(BC_n,B_n)$ \cite{Ko}
\begin{align}
P^{(BC_n,B_n)}_\lambda(x|a,b;q,t)=P_\lambda(x|-b^{1/2},ab^{1/2},-q^{1/2},q^{1/2};q,t).
\end{align} 
Setting $b=1$, we have the Macdonald polynomial of type $B_n$
\begin{align}
P^{(B_n)}_\lambda(x|a;q,t)=P^{(BC_n,B_n)}_\lambda(x|a,1;q,t).
\end{align}
Setting further $a=1$, we have  the Macdonald polynomial of type $D_n$
\begin{align}
P^{(D_n)}_\lambda(x|q,t)=P^{(BC_n,B_n)}_\lambda(x|1,1;q,t).
\end{align}

The application of the twist $\widehat{*}$ on the parameters 
$ (-b^{1/2},ab^{1/2},-q^{1/2},q^{1/2};q,t)$ gives
\begin{eqnarray}
(-(qb/t)^{1/2},a (qb/t)^{1/2},-q t^{-1/2},q t^{-1/2};q,q/t).\label{param-BD}
\end{eqnarray}

\subsection{Lassalle's conjectures}\label{Lassalle-con}
For the readers' convenience, we recall Lassalle's conjectures for Macdonald polynomials of type $B$, $C$ and $D$ 
with one row diagram.

Set 
\begin{align}
&
g^{(C_n)}_r(x)={(t;q)_r\over (q;q)_r}P^{(C_n)}_{(r)}(x|b;q,t)=
{(t;q)_r\over (q;q)_r}
P_{(r)}(x|-b^{1/2},b^{1/2}, -q^{1/2}b^{1/2},q^{1/2}b^{1/2}|q,t),\\
&
g^{(B_n)}_r(x)={(t;q)_r\over (q;q)_r}P^{(B_n)}_{(r)}(x|a;q,t)=
{(t;q)_r\over (q;q)_r}
P_{(r)}(x|-1,a, -q^{1/2},q^{1/2}|q,t),\\
&
g^{(D_n)}_r(x)={(t;q)_r\over (q;q)_r}P^{(D_n)}_{(r)}(x|q,t)=
{(t;q)_r\over (q;q)_r}
P_{(r)}(x|-1,1, -q^{1/2},q^{1/2}|q,t).
\end{align}
Lassalle's conjectures \cite[p.8, Conjecture 1, p.10, Conjecture 3, p.11, Conjecture 4]{Lass}  read
\begin{align}
&
g_r^{(C_n)}= 
\sum_{i=0}^{\lfloor r/2 \rfloor}
G_{r-2i} t^i 
{(b/t;q)_i\over (q;q)_i} 
{(t^n q^{r-i};q)_i \over (b t^{n-1}q^{r-i};q)_i} {1-t^n q^{r-2i}\over 1-t^n q^{r-i}} ,\\
&
g_r^{(D_n)}= 
\sum_{i=0}^{\lfloor r/2 \rfloor}
G_{r-2i} t^i 
{(1/t;q)_i\over (q;q)_i} 
{(t^n q^{r-i};q)_i \over ( t^{n-1}q^{r-i};q)_i} {1-t^n q^{r-2i}\over 1-t^n q^{r-i}} ,\\
&g_r^{(B_n)}= 
\sum_{i=0}^{r}
g_{r-i}^{(D_n)} 
{(a;q)_i\over (q;q)_i} 
{(t^n q^{r-i};q)_i \over (t^{n-1}q^{r-i};q)_i} 
{(t^{2n-2}q^{2r-i+1};q)_i \over 
(a t^{2n-2}q^{2r-i};q)_i}.
\end{align}
Namely, for type $C$ and $D$ we have
\begin{align}
&
g_r^{(C_n)}= 
\sum_{i=0}^{\lfloor r/2 \rfloor}
G_{r-2i} 
{(b/t;q)_i\over (q;q)_i} 
{(t^{-n} q^{-r};q)_i \over (t^{-n+1}q^{-r+1}/b;q)_i} {1-t^{-n} q^{-r+2i}\over 1-t^{-n} q^{r}} (t^2/qb)^i ,\\
&
g_r^{(D_n)}= 
\sum_{i=0}^{\lfloor r/2 \rfloor}
G_{r-2i} 
{(1/t;q)_i\over (q;q)_i} 
{(t^{-n} q^{-r};q)_i \over (t^{-n+1}q^{-r+1};q)_i} {1-t^{-n} q^{-r+2i}\over 1-t^{-n} q^{r}} (t^2/q)^i .
\end{align}
For type $B$ we have
\begin{align}
&g_r^{(B_n)}= 
\sum_{i=0}^{r}
g_{r-i}^{(D_n)} 
{(a;q)_i\over (q;q)_i} 
{(t^{-n} q^{-r+1};q)_i \over (t^{-n+1}q^{-r};q)_i} 
{(t^{-2n+2}q^{-2r};q)_i \over 
(t^{-2n+2}q^{-2r+1}/a;q)_i}(t/a)^i\nonumber \\
&=
\sum_{i=0}^{r}
\sum_{j=0}^{\lfloor (r-i)/2 \rfloor}
G_{r-i-2j} 
{(a;q)_i\over (q;q)_i} 
{(t^{-n} q^{-r+1};q)_i \over (t^{-n+1}q^{-r};q)_i} 
{(t^{-2n+2}q^{-2r};q)_i \over 
(t^{-2n+2}q^{-2r+1}/a;q)_i}(t/a)^i\\
&\quad \times 
{(1/t;q)_j\over (q;q)_j} 
{(t^{-n} q^{-r+i};q)_j \over (t^{-n+1}q^{-r+i+1};q)_j} {1-t^{-n} q^{-r+i+2j}\over 1-t^{-n} q^{r+i}} (t^2/q)^j \nonumber \\
&=
\sum_{i=0}^{r}
\sum_{j=0}^{\lfloor (r-i)/2 \rfloor}
G_{r-i-2j} 
{(a;q)_i\over (q;q)_i} 
{(t^{-n+1} q^{-r+1};q)_i \over (t^{-n+1}q^{-r};q)_i} 
{(t^{-2n+2}q^{-2r};q)_i \over 
(t^{-2n+2}q^{-2r+1}/a;q)_i}(t/a)^i \nonumber\\
&\quad \times 
{(1/t;q)_j\over (q;q)_j} 
{(t^{-n} q^{-r};q)_{i+j} \over (t^{-n+1}q^{-r+1};q)_{i+j}} {1-t^{-n} q^{-r+i+2j}\over 1-t^{-n} q^{r}} (t^2/q)^j . \nonumber
\end{align}

\subsection{Conjecture about the Macdonald polynomial of type $B_2$}\label{con-B2}
We present a conjecture about the formal series $f(x;s)$ for type $B_2$. 
Let $\varepsilon_1,\varepsilon_2$ be the standard basis for ${\bf R}^2$.
The simple roots (for $B_2$) are 
$\alpha_1=\varepsilon_1-\varepsilon_2$,
$\alpha_2=\varepsilon_2$, fundamental weights are 
$\omega_1=\varepsilon_1$, 
$\omega_2=(\varepsilon_1+\varepsilon_2)/2$.
Let $P$ and $P^+$ be the weight lattice and the
cone of dominant weights.
Let $W$ be the Weyl group of type $B_2$. Denote the monomial
symmetric polynomials by $m_\lambda$ ($\lambda\in P^+$):
$
m_\lambda=\sum_{\mu\in W^\lambda} e^\mu.
$
We write $x_1=e^{\varepsilon_1}$, $x_2=e^{\varepsilon_2}$ for simplicity.  

Let $q$, $t$ and $T$ be indeterminates.
The Macdonald difference operator fortype  $B_2$ is defined by
\begin{align}
E_{\omega_1}(q,t,T)
&=
{1-t x_1/x_2\over 1-x_1/x_2}
{1-t  x_1x_2\over 1-x_1x_2}
{1-T x_1\over 1-x_1}
T_{q,x_1}+
{1-t x_2/x_1\over 1-x_2/x_1}
{1-t x_1x_2\over 1-x_1x_2}
{1-T x_2\over 1-x_2}
T_{q,x_2}\\
&+
{1-t /x_1x_2\over 1-1/x_1x_2}
{1-t  x_2/x_1\over 1-x_2/x_1}
{1-T /x_1\over 1-1/x_1}
T_{q^{-1},x_1}
+
 {1-t /x_1x_2\over 1-1/x_1x_2}
{1-t x_1/x_2\over 1-x_1/x_2}
{1-T /x_2\over 1-1/x_2}
T_{q^{-1},x_2}.\nonumber
\end{align}
The Macdonald polynomials $P_\lambda(x_1,x_2;q,t,T)$ of type $B_2$ 
are uniquely characterized by the 
following conditions.
\begin{align}
(i)& \,\,P_\lambda(x_1,x_2;q,t,T)=m_\lambda+\sum_{\mu\in P^+,\mu<\lambda}
a_{\lambda\mu}(q,t,T) m_\mu,\\
(ii)& \,\,E_{\omega_1}(q,t,T)P_\lambda(x_1,x_2;q,t,T)=
c_\lambda P_\lambda(x_1,x_2;q,t,T),
\end{align}
where 
\begin{align}
c_{r_1\omega_1+r_2\omega_2}
=t^2T q^{r_1+r_2/2}+
t T q^{r_2/2}+t q^{-r_2/2}+
q^{-r_1-r_2/2}\qquad (r_1,r_2\in \bbZ_{\geq 0}),
\end{align}
or equivalently
\begin{align}
c_{\lambda_1\varepsilon_1+\lambda_2\varepsilon_2}
=t^2T q^{\lambda_1}+
t T q^{\lambda_2}+t q^{-\lambda_2}+
q^{-\lambda_1}\quad (\lambda_1,\lambda_2\in\bbZ_{\geq 0}/2,\lambda_1+\lambda_2\in \bbZ_{\geq 0},
\lambda_1\geq \lambda_2).
\end{align}

Let $s_1,s_2$ be generic parameters. 
Introduce variables $\lambda_1,\lambda_2$ satisfying $s_1=t T^{1/2} q^{\lambda_1},s_2=T^{1/2} q^{\lambda_2}$.
Note that we have $T_{q,x_1} x_1^{\lambda_1}= t^{-1}T^{-1/2} s_1  x_1^{\lambda_1},
T_{q,x_2} x_2^{\lambda_2}= T^{-1/2} s_1  x_1^{\lambda_1}$.

Set 
\begin{align}
&f^{B_2}(x_1,x_2;s_1,s_2,q,t,T)
= x_1^{\lambda_1}x_2^{\lambda_2}\sum_{n=0}^\infty f^{B_2}_n(x_1,x_2;s_1,s_2,q,t,T),\label{f_n_B}\\
&f^{B_2}_n(x_1,x_2;s_1,s_2,q,t,T)\label{f-B2}\\
&=
(q^2 t/T^2)^n 
{(q/t;q)_n (T/t;q)_n (T;q)_n (T/s_1s_2;q)_n (T/t s_1 s_2;q)_n (q/t s_1 s_2;q)_n
\over 
(q;q)_n (q/s_1^2;q)_n (q/s_2^2;q)_n (qs_1/s_2;q)_n (q s_2/s_1;q)_n
(q/s_1 s_2;q)_n}\nonumber \\
&
\times {
(T/s_1^2;q)_{2n} (T/s_2^2;q)_{2n}
\over 
(T/s_1s_2;q)_{2n}(T/t s_1s_2;q)_{2n}} 
(1/x_1 x_2)^n\nonumber \\
&\times
\sum_{\theta_{1},\theta_{2},\theta_{3},\theta_4\geq 0}
c^{B_2}(n,\theta_{1},\theta_{2}, \theta_{3},\theta_4;s_1,s_2,q,t,T)
(x_2/x_1)^{\theta_{1}}(1/x_2)^{\theta_{2}}
(1/x_1)^{\theta_{3}}(1/x_1x_2)^{\theta_{4}},\nonumber
\end{align}
where
\begin{align}
&
c^{B_2}(n,\theta_{1},\theta_{2}, \theta_{3},\theta_4;s_1,s_2,q,t,T)\nonumber\\
&=
(q/t)^{\theta_{1}} 
{(t;q)_{\theta_{1}} \over (q;q)_{\theta_{1}} }
{(q^{\theta_{3}-\theta_{2}}t s_2/s_1;q)_{\theta_{1}} \over 
(q^{\theta_{3}-\theta_{2}}q s_2/s_1;q)_{\theta_{1}} } \label{c-B2}\\
&\times
(q/T)^{\theta_{2}} 
{(q^n T;q)_{\theta_{2}} \over (q;q)_{\theta_{2}} }
{(q^{2n} T/s_2^2;q)_{\theta_{2}} \over 
(q^{n}q/s_2^2;q)_{\theta_{2}} }
{(q^n T/s_1s_2;q)_{\theta_2} \over (q^{2n} T/s_1s_2;q)_{\theta_2}}
{(q s_1/s_2;q)_{\theta_2} \over (q^{n} q s_1/s_2;q)_{\theta_2}} \nonumber\\
&\times
(q/T)^{\theta_{3}} 
{(q^n T;q)_{\theta_{3}} \over (q;q)_{\theta_{3}} }
{(q^{2n} T/s_1^2;q)_{\theta_{3}} \over 
(q^{n}q/s_1^2;q)_{\theta_{3}} } 
{(t s_2/s_1;q)_{\theta_{3}} \over 
(q^n qs_2/s_1;q)_{\theta_{3}} }
{(q^{-\theta_{2}}q s_2/t s_1;q)_{\theta_{3}} \over 
(q^{-\theta_{2}}s_2/s_1;q)_{\theta_{3}} }\nonumber\\
&\times
{(q^n T/ s_1s_2;q)_{\theta_{3}} \over 
(q^{2n} qT/t s_1 s_2;q)_{\theta_{3}} }
{(q^{2n} q^{\theta_{2}}q T/t s_1s_2;q)_{\theta_{3}} \over 
(q^{2n} q^{\theta_{2}}T/s_1 s_2;q)_{\theta_{3}} }\nonumber\\
&\times
(q/t)^{\theta_{4}} 
{(t;q)_{\theta_{4}} \over (q;q)_{\theta_{4}} }
{(q^{2n}q^{\theta_{2}+\theta_{3}} tT/ s_1s_2;q)_{\theta_{4}} \over 
(q^{2n}q^{\theta_{2}+\theta_{3}} qT/ s_1s_2;q)_{\theta_{4}} }.\nonumber
\end{align}

\begin{con}\label{conB}
The series $f^{B_2}(x_1,x_2;s_1,s_2,q,t,T)$ in (\ref{f_n_B}) satisfies
\begin{align}
E_{\omega_1}(q,t,T)f^{B_2}(x_1,x_2;s_1,s_2,q,t,T)
=t T^{1/2}
(s_1+s_2+s_1^{-1}+s_2^{-1})
f^{B_2}(x_1,x_2;s_1,s_2,q,t,T). \label{difeq-B}
\end{align}
\end{con}

\begin{con}\label{conB2}
Let $r_1,r_2 \in \bbZ_{\geq 0}$. We have
\begin{eqnarray}
P_{r_1\omega_1+r_2\omega_2}(x_1,x_2;q,t,T)
=x_1^{r_1+r_2/2}x_2^{r_2/2}
 f^{B_2}(x_1,x_2;t T^{1/2}q^{r_1+r_2/2},T^{1/2} q^{r_2/2},q,t,T).
\end{eqnarray}
Or equivalently, for any half-partition $(\lambda_1,\lambda_2)$, we have
\begin{eqnarray}
P_{\lambda_1\varepsilon_1+\lambda_2\varepsilon_2}(x_1,x_2;q,t,T)
=x_1^{\lambda_1}x_2^{\lambda_2}
 f^{B_2}(x_1,x_2;t T^{1/2}q^{\lambda_1},T^{1/2} q^{\lambda_2},q,t,T).
\end{eqnarray}
\end{con}

\begin{rmk}
1) When the parameters are specialized as 
 \begin{align}
 s_1=t T^{1/2}q^{r_1+r_2/2},
 \qquad s_2=T^{1/2} q^{r_2/2} \qquad (r_1,r_2\in \bbZ_{\geq 0}). \label{special-s}
 \end{align}
the series $f^{B_2}(x_1,x_2;s_1,s_2,q,t,T)$ 
becomes  truncated.  It also has to be symmetric with respect to 
the action of the Weyl group of type $B_2$. The truncation can be checked explicitly from (\ref{c-B2}).
We have not checked the symmetry yet.\\
  2) If $r=r \omega_1$ ($s_1=tT^{1/2} q^r,s_2=T^{1/2}$, $r\in \bbZ_{\geq 0},r_2=0$),
 in view of (\ref{f-B2}),
 we have $f^{B_2}_n(x_1,x_2;s_1,s_2,q,t,T)=0$ when $n>0$. 
 Hence Conjecture \ref{conB2} implies the  threefold summation formula
 \begin{align}
 &P_{r\omega_1}(x_1,x_2;q,t,T)=
  x_1^{r}
 f^{B_2}_0(x_1,x_2;t T^{1/2}q^{r},T^{1/2},q,t,T)\nonumber\\
&
=
 x_1^{r}
\sum_{\theta_{1},\theta_{3},\theta_4\geq 0}
(x_2/x_1)^{\theta_{1}}
(1/x_1)^{\theta_{3}}(1/x_1x_2)^{\theta_{4}}\nonumber\\
&\times
(q/t)^{\theta_{1}} 
{(t;q)_{\theta_{1}} \over (q;q)_{\theta_{1}} }
{(q^{\theta_{3}-r};q)_{\theta_{1}} \over 
(q^{\theta_{3}-r+1}/t;q)_{\theta_{1}} }
(q/T)^{\theta_{3}} 
{(T;q)_{\theta_{3}} \over (q;q)_{\theta_{3}} }
{(q^{-2r}/t^2;q)_{\theta_{3}} \over 
(q^{-2r+1}/t^2T;q)_{\theta_{3}} }\\
&\times
{(q^{-r};q)_{\theta_{3}} \over 
(q^{-r+1}/t;q)_{\theta_{3}} }
{(q^{-r+1}/t^2;q)_{\theta_{3}} \over 
(q^{-r}/t;q)_{\theta_{3}} }
(q/t)^{\theta_{4}} 
{(t;q)_{\theta_{4}} \over (q;q)_{\theta_{4}} }
{(q^{\theta_{3}-r} ;q)_{\theta_{4}} \over 
(q^{\theta_{3}-r+1}/t;q)_{\theta_{4}} }.\nonumber
\end{align}
 3)
Let $\lambda=r_1\omega_1+r_2\omega_2$ ($r_1,r_2\in \bbZ_{\geq 0}$), and 
$s_1,s_2$ are specialized as in  (\ref{special-s}).
If we further specialize the parameters as $q=t=T$, 
the irreducible character 
${\rm ch}_{r_1\omega_1+r_2\omega_2}$
for the Lie algebra of type  
$B_2$
must be recovered from Conjecture  \ref{conB2}.
By checking termination of the series, one finds that Conjecture \ref{conB2} implies 
\begin{eqnarray}
{\rm ch}_{r_1\omega_1+r_2\omega_2}=
x_1^{r_1+r_2/2}x_2^{r_2/2}
\sum_{\theta_1,\theta_2,\theta_3,\theta_4\geq 0, \atop
{\theta_1-\theta_2+\theta_3\leq r_1, \atop{
\theta_2\leq r_2, \atop {
\theta_3\leq r_1, \atop
\theta_2+\theta_3+\theta_4\leq r_1+r_2}}}}
(x_2/x_1)^{\theta_{1}}(1/x_2)^{\theta_{2}}
(1/x_1)^{\theta_{3}}(1/x_1x_2)^{\theta_{4}}.
\end{eqnarray}
4) Conjecture \ref{conB2} has been checked up to 
 $r_1+r_2\leq 6$.

\end{rmk}

\end{document}